\documentclass[review]{siamart190516}   % Use [final] when ready.

\usepackage{pdflscape}
\usepackage{siampaper}

\newcommand{\USYMLQR}{\textsc{Usymlqr}\xspace}

\newcommand{\SYMMLQ}{\textsc{Symmlq}\xspace}
\newcommand{\CGM}{\textsc{Cg}\xspace}

\newcommand{\MINRES}{\textsc{Minres}\xspace}

\newcommand{\BLCG}{\textsc{Block-Cg}\xspace}
\newcommand{\BLMINRES}{\textsc{Block-Minres}\xspace}

\newcommand{\MINRESQLP}{\textsc{Minres-qlp}\xspace}

\newcommand{\LSLQ}{\textsc{Lslq}\xspace}
\newcommand{\LSQR}{\textsc{Lsqr}\xspace}
\newcommand{\LSMR}{\textsc{Lsmr}\xspace}

\newcommand{\LNLQ}{\textsc{Lnlq}\xspace}
\newcommand{\CRAIG}{\textsc{Craig}\xspace}
\newcommand{\CRAIGMR}{\textsc{Craigmr}\xspace}

\newcommand{\TriCG}{\textsc{TriCG}\xspace}
\newcommand{\TriMR}{\textsc{TriMR}\xspace}
\newcommand{\TriLQ}{\textsc{TriLQ}\xspace}

\def\T{^T\!}

\newlength{\forwidth}
\newcommand{\bmat}[1]{\begin{bmatrix} #1 \end{bmatrix}}

\usepackage[utf8]{inputenc}
\usepackage[T1]{fontenc}

\graphicspath{{graphics/}}  % path to figures
\usepackage{tikzscale}
\usepackage{pgfplots}
\pgfplotsset{compat=newest}
\usetikzlibrary{external}
\usetikzlibrary{math}

\newcommand*{\includetikzgraphics}[2][]{%
  \tikzsetnextfilename{#2}%
  \includegraphics[#1]{#2.tikz}
}

\makeatletter
\renewcommand{\todo}[2][]{\tikzexternaldisable\@todo[#1]{#2}\tikzexternalenable}
\makeatother

\newcommand{\cahiernumber}{41}  % Insert your Cahier du GERAD number.

\usepackage{kbordermatrix}

\newcommand{\diag}{\mathop{\text{diag}}}
\title{%
  {\TriCG and \TriMR}:\@ Two Iterative Methods for Symmetric Quasi-Definite Systems
}

\author{%
  Alexis Montoison%
  \thanks{%
    GERAD and Department of Mathematics and Industrial Engineering,
    Polytechnique Montr\'eal, QC, Canada.
    E-mail: \mailto{alexis.montoison@polymtl.ca}.
    Research supported by an excellence scholarship of the IVADO institute.
  }
  \and
  Dominique Orban%
  \thanks{%
    GERAD and Department of Mathematics and Industrial Engineering,
    Polytechnique Montr\'eal, QC, Canada.
    E-mail: \mailto{dominique.orban@gerad.ca}.
    Research partially supported by an NSERC Discovery Grant.
  }
}
\date{\today}

\begin{document}

  \nolinenumbers
  \maketitle

  \thispagestyle{firstpage}
  \pagestyle{myheadings}

  \begin{abstract}
    We introduce iterative methods named \TriCG and \TriMR for solving symmetric quasi-definite systems based on the orthogonal tridiagonalization process proposed by Saunders, Simon and Yip in 1988.
    % and its equivalent block-Lanczos formulation.
    \TriCG and \TriMR are tantamount to preconditioned \BLCG and \BLMINRES with two right-hand sides in which the two approximate solutions are summed at each iteration, but require less storage and work per iteration.
    % We state optimality properties of the iterates of each method.
    We evaluate the performance of \TriCG and \TriMR on linear systems generated from the SuiteSparse Matrix Collection and from discretized and stablized Stokes equations.
    We compare \TriCG and \TriMR with \SYMMLQ and \MINRES, the recommended Krylov methods for symmetric and indefinite systems.
    In all our experiments, \TriCG and \TriMR terminate earlier than \SYMMLQ and \MINRES on a residual-based stopping condition with an improvement of up to 50\% in terms of number of iterations.
    They also terminate more reliably than \BLCG and \BLMINRES.
    Experiments in quadruple and octuple precision suggest that loss of orthogonality in the basis vectors is significantly less pronounced in \TriCG and \TriMR than in \BLCG and \BLMINRES.
  \end{abstract}


  % https://mathscinet.ams.org/mathscinet/msc/pdfs/classifications2020.pdf

  %\tableofcontents
  %\listoftodos\relax

  \section{Introduction}

  Consider a linear system of the form
  \begin{equation}
    \label{eq:sqd}
    \begin{bmatrix} M & \phantom{-}A \\ A^T & -N \end{bmatrix}
    \begin{bmatrix} x \\ y \end{bmatrix}
    =
    \begin{bmatrix} b \\ c \end{bmatrix},
  \end{equation}
  where $M \in \R^{m \times m}$ and $N \in \R^{n \times n}$ are symmetric positive definite, $b \in \R^m$ and $c \in \R^n$ are not simultaneously zero, and $A \in \R^{m \times n}$ can have any shape.

  We use the definition given by \cite{orban-arioli-2017}: a matrix $K$ is symmetric quasi-definite (SQD) if $K=K^T$ and there exists a permutation matrix $P$ such that $P\T K P$ has the form~\eqref{eq:sqd}.
  In particular,~\eqref{eq:sqd} arises in interior-point methods for inequality-constrained optimization \citep{altman-gondzio-1999, friedlander-orban-2012}
  and in the context of stabilized mixed finite elements methods \citep{elman-silvester-wathen-2014}.
  SQD matrices are indefinite and nonsingular \citep{vanderbei-1995}.

  In this paper, we develop two iterative methods named \TriCG and \TriMR specialized for~\eqref{eq:sqd}, based on the orthogonal tridiagonalization process in elliptic norms and a closely-related formulation as a preconditioned block-Lanczos method. % block-\CGM \citep{leary-1980}.
  Relations between both processes are described in detail and we show how the preconditioned block-Lanczos process with two specific right-hand sides can generate structured Krylov bases with SQD systems.
  Our main motivation for developing \TriCG and \TriMR comes from an absence of iterative methods for~\eqref{eq:sqd} that exploit the SQD structure when both $b$ and $c$ are nonzero.

  The rest of this paper is organized as follows.
  We state the defining properties of \TriCG of \TriMR and describe their implementations in detail.
  In a second stage, we compare \TriCG and \TriMR with their block counterparts \BLCG and \BLMINRES as well as \SYMMLQ and \MINRES on two set of problems.
  The first set uses the SuiteSparse Matrix Collection without preconditioning.
  The second set is composed of discretized and stablized Stokes equations and requires preconditioning.
  Finally, we discuss extensions of \TriCG and \TriMR and their uses outside the context of SQD systems.

  \subsection*{Related research}

  \cite{orban-arioli-2017} expose the state of the art on iterative methods for SQD systems.
  They explain that existing Krylov methods for symmetric indefinite systems, such as \SYMMLQ and \MINRES \citep{paige-saunders-1975} or \MINRESQLP \citep{choi-paige-minres-2011} do not exploit the rich structure of~\eqref{eq:sqd}.
  \CGM \citep{hestenes-stiefel-1952} is well defined on~\eqref{eq:sqd} provided $b=0$ or $c=0$.
  Otherwise, it may break down.
  For example, any SQD system such that $b\T M b + 2 b\T A c - c\T N c = 0$ causes breakdown on the first \CGM iteration.
  \cite{orban-arioli-2017} show that preconditioned and regularized \LSQR \citep{paige-saunders-1982} and \LSMR \citep{fong-saunders-2011} take advantage of block structure to solve
  \begin{equation}
    \label{eq:sqd_b}
    \begin{bmatrix} M & \phantom{-}A \\ A^T & -N \end{bmatrix}
    \begin{bmatrix} x \\ y \end{bmatrix}
    =
    \begin{bmatrix} b \\ 0 \end{bmatrix},
  \end{equation}
  and are equivalent to \CGM and \MINRES applied to the normal equations
  \begin{equation}
    \label{eq:normal_equation1}
    (A\T M^{-1} A  + N) y = A\T M^{-1} b \quad \text{with} \quad x = M^{-1}(b - Ay).
  \end{equation}
  They also show that preconditioned and regularized \CRAIG \citep{craig-1955} and \CRAIGMR \citep{orban-arioli-2017} solve
  \begin{equation}
    \label{eq:sqd_c}
    \begin{bmatrix} M & \phantom{-}A \\ A^T & -N \end{bmatrix}
    \begin{bmatrix} x \\ y \end{bmatrix}
    =
    \begin{bmatrix} 0 \\ c \end{bmatrix},
  \end{equation}
  and are equivalent to \CGM and \MINRES applied to the Schur-complement equations
  \begin{equation}
    \label{eq:normal_equation2}
    (A\T M^{-1} A + N) y = -c \quad \text{with} \quad x = -M^{-1} A y.
  \end{equation}
  In a similar vein, \cite{estrin-orban-saunders-2019a,estrin-orban-saunders-2019b} develop \LSLQ and \LNLQ to solve~\eqref{eq:sqd_b} and~\eqref{eq:sqd_c}, respectively, and explain that preconditioned and regularized variants of those methods are equivalent to \SYMMLQ applied to~\eqref{eq:normal_equation1} and to~\eqref{eq:normal_equation2}.

  When \(b \neq 0\) and \(c \neq 0\), one possibility is to shift the right-hand side to recover~\eqref{eq:sqd_b} or~\eqref{eq:sqd_c}.
  For instance, we can solve $-N \Delta y = c$ and add $(0, \Delta y)$ to the solution of~\eqref{eq:sqd_b} with right-hand side $(b - A \Delta y, 0)$.
  \cite{buttari-orban-ruiz-titley_peloquin-2019} developed \USYMLQR upon the orthogonal tridiagonalization process \citep{saunders-simon-yip-1988}, henceforth the \emph{SSY process}, for the saddle-point system
  \begin{equation}
    \label{eq:saddle_point}
    \begin{bmatrix} M & A \\ A^T & 0 \end{bmatrix}
    \begin{bmatrix} x \\ y \end{bmatrix}
    =
    \begin{bmatrix} b \\ c \end{bmatrix}.
  \end{equation}

  The SSY process requires two initial vectors, which makes it ideal to develop \TriCG and \TriMR, two methods specialized for SQD systems in the case where both $b$ and $c$ are nonzero.
  Its close relation to a block-Lanczos process induces similarities between \TriCG and the block-\CGM method \citep{leary-1980} as well as \TriMR and block-\MINRES methods.
  The block-Lanczos process was initialy developed to compute eigenvalues \citep{golub-underwood-1977}.
  Thereafter it was used to find nullspaces \citep{montgomery-1995} and solve linear systems with multiple right-hand sides \citep{guennouni-jbilou-sadok-2004}.
  In this paper, a novel application of this process is presented for SQD linear systems with a single right-hand side.

  \subsection*{Notation}

  Vectors and scalars are denoted by lowercase Latin and Greek letters, respectively.
  Matrices are denoted by capital Latin letters, except for \(2\)\(\times\)\(2\) blocks, which are represented by capital Greek letters.
  Rare exceptions to those rules are indicated explicitly in the text.
  For a vector $v$, $\|v\|$ denotes the Euclidean norm of $v$, and for a symmetric and positive-definite matrix \(M\), the $M$-norm of $v$ is \(\|v\|_M^2 = v\T M v\).
  The shorthand $v \mapsto M \backslash v$ represents an operator that returns the solution of $Mu = v$.
  The vector $e_i$ is the $i$-th column of an identity matrix of size dictated by the context.
  \(I_k\) represents the \(k\)\(\times\)\(k\) identity operator.
  The shorthand $\diag(\omega_1, \cdots, \omega_k)$ represents the \(k\)\(\times\)\(k\) diagonal matrix with \(\omega_1, \dots, \omega_k\) on its diagonal.
  We denote by \(K\) the SQD matrix of~\eqref{eq:sqd} and
  \begin{equation}
    \label{eq:def-H-K0-B}
    K_0 :=
    \begin{bmatrix}
      0 & A \\
      A\T & 0
    \end{bmatrix},
    \qquad
    H := \blkdiag(M, N) =
    \begin{bmatrix}
      M & 0 \\
      0 & N
    \end{bmatrix},
    \qquad
    B :=
    \begin{bmatrix}
      b & 0 \\
      0 & c
    \end{bmatrix}.
  \end{equation}
  We abuse notation and write $(b, c)$ to represent the column vector $\begin{bmatrix} b\T & c\T \end{bmatrix}^T$.

  \section{Processes}

  In this section, we state the SSY process in elliptic norms, which is the foundation for the development of the methods \TriCG and \TriMR, its relation with the preconditioned block-Lanczos process, and how they can accommodate regularization.

  \subsection{The orthogonal tridiagonalization process in elliptic norms}

  The SSY process generalized by \cite{buttari-orban-ruiz-titley_peloquin-2019} in terms of elliptic norms defined by positive definite $M$ and $N$ generates sequences of vectors ${v_k}$ and ${u_k}$ such that $v_i\T M v_j = \delta_{ij}$ and $u_i\T N u_j = \delta_{ij}$ in exact arithmetic.
  The process is stated as \Cref{alg:ssy}, where we use the shorthand notation \(\beta_1 M v_1 = b\) to summarize the normalization operations
  \begin{enumerate}
    \item set \(\bar{v}_1 = b\);
    \item solve \(M v_1 = \bar{v}_1\);
    \item compute \(\beta_1 = (\bar{v}_1\T v_1)^{\frac12}\);
    \item normalize \(\bar{v}_1 \leftarrow \bar{v}_1 / \beta_1\) and \(v_1 \leftarrow v_1 / \beta_1\),
  \end{enumerate}
  and similarly for \(\gamma_1 N u_1 = c\) and subsequent normalization steps.
  When $M$ and $N$ are not the identity, the above normalization operations only require solves with $M$ and $N$.
  Occurrences of $M v_k$ and $N u_k$ in a right-hand side in \Cref{alg:ssy} simply mean that we substitute $\bar{v}_k$ and $\bar{u}_k$, respectively, so that \(M\) and \(N\) themselves are not needed.
  The process terminates if there exists an index \(k\) such that \(\beta_{k+1} = 0\) or \(\gamma_{k+1} = 0\).

  \begin{algorithm}[ht]
    \caption{%
      Orthogonal Tridiagonalization Process in Elliptic Norms
    }
    \label{alg:ssy}
    \begin{algorithmic}[1]
      \Require $A$, $b$, $c$, $v \mapsto M \backslash v$, $u \mapsto N \backslash u$
      \State $v_0 = 0$, $u_0 = 0$
      \State $\beta_1 M v_1 = b$, $\gamma_1 N u_1 = c$ \Comment{$(\beta_1, \, \gamma_1) > 0$ so that $\|v_1\|_M = \|u_1\|_N = 1$}
      \For{$k$ = 1,~2,~\(\dots\)}
        \State $q = A u_k - \gamma_k M v_{k-1}$, $\alpha_k = v_k^T q$
        \State $p = A^T v_k - \beta_k N u_{k-1}$
        \State $\beta_{k+1} M v_{k+1} = q - \alpha_k M v_k$ \Comment{$\beta_{k+1} > 0$ so that $\|v_{k+1}\|_M = 1$}
        \State $\gamma_{k+1} N u_{k+1} = p - \alpha_k N u_k$ \Comment{$\gamma_{k+1} > 0$ so that $\|u_{k+1}\|_N = 1$}
      \EndFor
    \end{algorithmic}
  \end{algorithm}

  We denote \(V_k = \begin{bmatrix} v_1 & \dots & v_k \end{bmatrix}\) and \(U_k = \begin{bmatrix} u_1 & \dots & u_k \end{bmatrix}\).
  After \(k\) iterations of \Cref{alg:ssy}, the situation may be summarized as
  \begin{subequations}
    \label{eq:triortho}
    \begin{align}
         A U_k & = M V_k T_k~ + \beta_{k+1} M v_{k+1} e_k^T = M V_{k+1} T_{k+1,k}
         \label{eq:triortho-V}
      \\ A^T V_k & = N U_k T_k^T + \gamma_{k+1} N u_{k+1} e_k^T = N U_{k+1} T_{k,k+1}^T
        \label{eq:triortho-U}\\
        V_k\T M V_k &= U_k\T N U_k = I_k,
        \label{eq:triortho-UV}
    \end{align}
  \end{subequations}
  where
  \begin{equation*}
    T_k =
    \begin{bmatrix}
      \alpha_1 & \gamma_2 &          & \\
      \beta_2  & \alpha_2 & \ddots   & \\
               & \ddots   & \ddots   & \gamma_k \\
               &          & \beta_k  & \alpha_k
    \end{bmatrix},
    \qquad
    T_{k,k+1} =
    \begin{bmatrix}
      T_{k} & \gamma_{k+1} e_k
    \end{bmatrix},
    \qquad
    T_{k+1,k} =
    \begin{bmatrix}
      T_{k} \\
      \beta_{k+1} e_{k}^T
    \end{bmatrix}.
  \end{equation*}

  Equations~\eqref{eq:triortho-V}--\eqref{eq:triortho-U} hold to within machine precision despite loss of orthogonality, but~\eqref{eq:triortho-UV} and $V_k\T A U_k = T_k$ hold only in exact arithmetic.

  \subsection{Relation with preconditioned block-Lanczos process}

  \cite{saunders-simon-yip-1988} note Beresford Parlett's observation that the subspaces generated by \Cref{alg:ssy} in the Euclidean norm can be viewed as the union of those generated by the block-Lanczos process applied to $A\T A$ and $A A\T$ with respective starting blocks $\begin{bmatrix} c & A\T b \end{bmatrix}$ and $\begin{bmatrix} b & Ac \end{bmatrix}$.
  \cite{golub-stoll-wathen-2008} pushed the observation further in terms of the block-Lanczos process applied to \(K_0\) in \eqref{eq:def-H-K0-B}.
  This section summarizes the latter observations and incorporates the preconditioner \(H\).

  Pasting~\eqref{eq:triortho} together results in
  \begin{equation}
    \label{eq:ssy-bl}
    \begin{bmatrix} 0 & A \\ A^T & 0 \end{bmatrix}
    \begin{bmatrix} V_k & 0 \\ 0 & U_k \end{bmatrix}
    =
    \begin{bmatrix} M & 0 \\ 0 & N \end{bmatrix}
    \begin{bmatrix} V_{k+1} & 0 \\ 0 & U_{k+1} \end{bmatrix}
    \begin{bmatrix} 0 & T_{k+1,k} \\ T_{k,k+1}^T & 0 \end{bmatrix},
  \end{equation}
  which resembles a Krylov process in which basis vectors have been permuted.
  Let
  $$P_k := \bmat{e_1 \!&\! e_{k+1} \!&\! \cdots \!&\! e_i \!&\! e_{k+i} \!&\! \cdots \!&\! e_k \!&\! e_{2k}} = \bmat{E_1 \!&\! \cdots \!&\! E_k} \quad \text{where} \quad
    E_k = \bmat{e_k \!&\! \\ \!&\! e_k},$$
  denote the permutation, introduced by \cite{paige-1974}, that restores the order in which \cref{alg:ssy} generates basis vectors:
  \begin{equation}
    \label{eq:form-Wk}
    W_k := \begin{bmatrix} V_k & 0 \\ 0 & U_k \end{bmatrix} P_k
    =
    \begin{bmatrix}
      w_1 & \cdots & w_k
    \end{bmatrix}
    \quad \text{where} \quad
    w_k =
    \begin{bmatrix}
      v_k & 0 \\ 0 & u_k
    \end{bmatrix}.
  \end{equation}

  Although $w_k$ is a matrix, we use a lowercase letter due to the close link with the vectors $v_k$ and $u_k$.
  The projection of $K_0$ in the Krylov subspace $\Span\{w_1,\cdots, w_k\}$ is also shuffled to symmetric block-tridiagonal form with blocks of size \(2\):
  \begin{align}
    K_0 W_k
    &=
    \begin{bmatrix} M & 0 \\ 0 & N \end{bmatrix}
    \begin{bmatrix} V_{k+1} & 0 \\ 0 & U_{k+1} \end{bmatrix} P_{k+1} P_{k+1}^T
    \begin{bmatrix} 0 & T_{k+1,k} \\ T_{k,k+1}^T & 0 \end{bmatrix} P_k
    \nonumber \\
    &=
    H W_{k+1} F_{k+1,k},
    \label{eq:bl-A-A'}
  \end{align}
  where
  $$F_{k+1,k}
  =
  \begin{bmatrix}
      \Omega_1 & \Psi_2   &        & \\
      \Psi_2^T & \Omega_2 & \ddots & \\
               & \ddots   & \ddots & \Psi_k \\
               &          & \ddots & \Omega_k \\
               &          &        & \Psi_{k+1}^T
  \end{bmatrix},
  \qquad
  \Omega_k = \begin{bmatrix} 0 & \alpha_k \\ \alpha_k & 0 \end{bmatrix},
    \qquad
  \Psi_k = \begin{bmatrix} 0 & \gamma_k \\ \beta_k & 0 \end{bmatrix}.$$
  The two relations at line 2 of \Cref{alg:ssy} can be rearranged as
  \begin{equation}
    \label{eq:bl-bc}
    \bmat{M & 0 \\ 0 & N} \bmat{v_1 & 0 \\ 0 & u_1} \bmat{\beta_1 & 0 \\ 0 & \gamma_1} = \bmat{b & 0 \\ 0 & c}\quad \Longleftrightarrow \quad H w_1 \Psi_1^T = B.
  \end{equation}
  The identities~\eqref{eq:bl-A-A'} and~\eqref{eq:bl-bc} characterize the preconditioned block-Lanczos process applied to \(K_0\) with preconditioner \(H\) and initial block \(B\).
  We summarize the process as \Cref{alg:block_lanczos} where all \(w_k \in \R^{(n + m) \times 2}\) and \(\Psi_k \in \R^{2 \times 2}\) are determined such that both $w_k\T H w_k = I_2$ and the equations on lines 2 and 5 are verified.

  \algsetblock[Name]{For}{EndFor}{}{.75em}
  \begin{algorithm}[ht]
    \caption{%
      Preconditioned Block-Lanczos Process
    }
    \label{alg:block_lanczos}
    \begin{algorithmic}[1]
      \Require $K_0$, $B$, $w \mapsto H \backslash w$
      \State $w_0 = 0$
      \State $H w_1 \Psi_1^T = B$\label{alg:init}
      \For{$k$ = 1,~2,~\(\dots\)}
        \State $\Omega_k = w_k\T K_0 w_k$
        \State $H w_{k+1} \Psi_{k+1}^T = K_0 w_k - H w_k \Omega_k - H w_{k-1} \Psi_k$\label{alg:recurrence}
      \EndFor
    \end{algorithmic}
  \end{algorithm}
  Note that \Cref{alg:ssy} and \Cref{alg:block_lanczos} require operators that return the solution of systems with coefficient $M$, $N$ and $H$.
  A specificity of \Cref{alg:block_lanczos} is that $w_k$ and $\Psi_k$ are not unique. They are commonly determined from the Gram-Schmidt process: $w_k \Psi_k^T$ is the QR decomposition of the right-hand side on lines \ref{alg:init} and \ref{alg:recurrence} of \Cref{alg:block_lanczos}.
  For instance, $F_{k+1,k}$ is pentadiagonal when we force $\Psi_k = \diag(\beta_k, \gamma_k)$ for all \(k\), in which case the structure of $w_k$ is
  \begin{equation*}
    w_{k} = \begin{bmatrix} v_{k} & 0 \\ 0 & u_{k} \end{bmatrix}
    \quad \text{($k$ odd)} \quad \text{and} \quad
    w_{k} = \begin{bmatrix} 0 & v_{k} \\ u_{k} & 0 \end{bmatrix}
    \quad \text{($k$ even)}.
  \end{equation*}

  \subsection{Regularization of the preconditioned block-Lanczos process}

  \begin{theorem}
    \label{theorem:sparcity-block-lanczos}
    Given the SQD matrix \(K\) and block right-hand side \(B\), the preconditioned Krylov basis $W_k$ generated by \Cref{alg:block_lanczos} has the form~\eqref{eq:form-Wk} where the vectors \(u_k\) and \(v_k\) are the same as those generated by \Cref{alg:ssy} with initial vectors \(b\) and \(c\).
    In addition,
    \begin{equation}
      \label{eq:sqd-block-lanczos}
      K
      W_k
      =
      H
      W_{k+1} S_{k+1, k},
      \qquad
      S_{k+1,k} :=
      \begin{bmatrix}
        \Theta_1 & \Psi_2   &        & \\
        \Psi_2^T & \Theta_2 & \ddots & \\
                 & \ddots   & \ddots & \Psi_k \\
                 &          & \ddots & \Theta_k \\
                 &          &        & \Psi_{k+1}^T
      \end{bmatrix},
    \end{equation}
    where
    \[
      \Theta_k =
      \begin{bmatrix}
        1 & \alpha_k \\
        \alpha_k & -1
      \end{bmatrix}
      \quad \text{and} \quad
      \Psi_k =
      \begin{bmatrix}
        0 & \gamma_k \\
        \beta_k & 0
      \end{bmatrix}.
    \]
    The scalars $\alpha_k$, $\beta_k$ and $\gamma_k$ are those generated by \Cref{alg:ssy} when it is applied to $A$ with initial vectors $b$ and $c$.
  \end{theorem}

  \begin{proof}
    Observe that \(K = K_0 + \blkdiag(M, -N)\).
    \Cref{alg:block_lanczos} applied to $K_0$ generates sparse pairs $w_k$ as in~\eqref{eq:form-Wk} because of the equivalence with \Cref{alg:ssy}.
    The term \(\blkdiag(M, -N)\) can be seen as a regularization term:
    \begin{equation}
      \label{eq:regularization_MN}
      \begin{bmatrix}
        M & \phantom{-}0 \\
        0 & -N
      \end{bmatrix}
      w_k
      =
      H w_k \Lambda_k
      \quad \text{with} \quad
      \Lambda_k :=
      \begin{bmatrix}
        1 & \phantom{-}0 \\
        0 & -1
      \end{bmatrix}.
    \end{equation}
    The identities~\eqref{eq:bl-A-A'} and~\eqref{eq:regularization_MN} allow us to write
    \begin{equation}
      \label{eq:regularization_sqd}
      K
      W_k
      =
      H
      \left(
      W_k
      \begin{bmatrix}
        \Omega_1 + \Lambda_1 & \Psi_2   &          & \\
        \Psi_2^T             & \ddots   & \ddots   & \\
                             & \ddots   & \ddots   & \Psi_k \\
                             &          & \Psi_k^T & \Omega_k + \Lambda_k \\
      \end{bmatrix}
      + w_{k+1} \Psi_{k+1}^T
      \right),
    \end{equation}
    which amounts to~\eqref{eq:sqd-block-lanczos} because $\Theta_k = \Omega_k + \Lambda_k$.
    The Krylov basis $W_k$ is not modified; only the projection of $K$ in the Krylov subspace is updated.
  \end{proof}

  Because of \cref{theorem:sparcity-block-lanczos}, the Krylov basis $W_k$ generated by \Cref{alg:block_lanczos} must have the sparsity structure~\eqref{eq:form-Wk}, so that only \(u_k\) and \(v_k\) need be generated, and they may be generated directly from \Cref{alg:ssy}.
  In addition, products with $M$ and $N$ are not required to generate $W_k$, so that the computational cost per iteration is reduced and less storage is required compared to \Cref{alg:block_lanczos}.

  \section{Methods}

  In this section, we develop two methods based upon \Cref{alg:ssy} in which iterates have the form
  \begin{equation}
    \label{eq:xy}
    \begin{bmatrix} x_k \\ y_k \end{bmatrix} = W_k z_k,
  \end{equation}
  where $z_k \in \R^{2k}$ is defined by certain optimality properties.
  Thanks to ~\eqref{eq:bl-bc} and ~\eqref{eq:sqd-block-lanczos}, the residual of~\eqref{eq:sqd} at any iterate of the form~\eqref{eq:xy} can be written
  \begin{align}
    r_k & =
    \begin{bmatrix}
      b \\ c
    \end{bmatrix}
    -
    \begin{bmatrix}
      M & \phantom{-}A \\
      A^T & -N
    \end{bmatrix}
    \begin{bmatrix}
      x_k \\ y_k
    \end{bmatrix}
    \nonumber \\
    & =
    H \left( w_1 \bmat{\beta_1 \\ \gamma_1} - W_{k+1} S_{k+1,k} z_k \right)
    \nonumber \\
    & = H W_{k+1} ( \beta_1 e_1 + \gamma_1 e_2 - S_{k+1,k} z_k ).
    \label{eq:residual_general}
  \end{align}
  In the next few sections, the particular choice of \(z_k\) yields a simplified expression for the residual.

  \subsection{Derivation of \TriCG}

  The \(k\)-th \TriCG iterate has the form~\eqref{eq:xy} with \(z_k\) defined by the Galerkin condition
  \begin{equation}
    \label{eq:galerkin-tricg}
    W_k^T r_k = 0
    \quad \Longleftrightarrow \quad
    W_k^T \left(\begin{bmatrix} b \\ c \end{bmatrix}
    -
    \begin{bmatrix} M & \phantom{-}A \\ A^T & -N \end{bmatrix}
    \begin{bmatrix} x_k \\ y_k \end{bmatrix}\right) = 0,
  \end{equation}
  which, thanks to~\eqref{eq:residual_general}, can be written as
  \[
    W_k^T H W_{k+1} \left(\beta_1 e_1 + \gamma_1 e_2 - S_{k+1,k} z_k \right) = 0.
  \]
  By construction of the Krylov basis, $W_k^T H W_k = I_{2k}$ and $w_i^T H w_j = 0$ for $i \neq j$ in exact arithmetic.
  Let $S_k \in \R^{2k \times 2k}$ denote the leading \((2k)\)\(\times\)\((2k)\) submatrix of $S_{k+1,k}$.
  This gives the \TriCG subproblem:
  \begin{equation}
    \label{eq:sub-tricg}
    S_k z_k = \beta_1 e_1 + \gamma_1 e_2.
  \end{equation}

  \subsubsection{Relation between \TriCG and block-\CGM}

  The \(k\)-th block-\CGM iterate is defined by the block-Galerkin condition
  \begin{equation}
    \label{eq:galerkin-block-cg}
    W_k^T \left(\begin{bmatrix} b & 0 \\ 0 & c \end{bmatrix}
    -
    \begin{bmatrix} M & \phantom{-}A \\ A^T & -N \end{bmatrix}
    \begin{bmatrix} x_k^b & x_k^c \\ y_k^b & y_k^c \end{bmatrix}\right) = 0,
  \end{equation}
  where $(x_k^b,~y_k^b) = W_k z_k^b$ and $(x_k^c,~y_k^c) = W_k z_k^c$.
  Accordingly, the \(k\)-th block-\CGM subproblem is
  \begin{equation}
  \label{eq:sub-block-cg}
  S_k \begin{bmatrix} z_k^b & z_k^c \end{bmatrix} = \begin{bmatrix} \beta_1 e_1 & \gamma_1 e_2 \end{bmatrix},
  \end{equation}
  so that $z_k^b$ and $z_k^c$ solve the subproblem associated with right-hand sides $(b,0)$ and $(0,c)$.
  The solutions of~\eqref{eq:sub-tricg} and~\eqref{eq:sub-block-cg} are connected via $z_k = z_k^b + z_k^c$, and the \TriCG and block-\CGM approximations are connected via $x_k = x_k^b + x_k^c$ and $y_k = y_k^b + y_k^c$.

  \subsubsection{An \texorpdfstring{$\mathrm{\mathbf{LDL^T}}$}{LDL} factorization}

  The connection between \Cref{alg:ssy} and \Cref{alg:block_lanczos} induces
  \begin{equation}
  \label{eq:Sk-Pk}
  S_k = P_k \begin{bmatrix} I_k & \phantom{-}T_k \\ T_k^T & -I_k \end{bmatrix} P_k^T,
  \end{equation}
  so that $S_k$ is SQD,
  and therefore nonsingular, and~\eqref{eq:sub-tricg} has a unique solution.
  Contrary to standard \CGM, the \TriCG iterates are always well-defined.
  \cite{vanderbei-1995} proved that SQD matrices are strongly factorizable, which means that,
  in particular, the factorization $S_k = L_k D_k L_k^T$ where $L_k$ is unit lower triangular and $D_k$ is diagonal always exists.
  Subsequently, the solution $z_k$ of~\eqref{eq:sub-tricg} can be determined via forward and backward sweeps, although the next section shows that computing \(z_k\) is not necessary.
  The factorization of $S_k$ can be updated at each iteration.
  Let
  \[
    D_k = \begin{bmatrix}
            d_1 &        &       \\
                & \ddots &       \\
                &        & d_{2k}
          \end{bmatrix}\!,
    \mathhfill
    L_k = \begin{bmatrix}
            \Delta_1 &          &          &          \\
            \Gamma_2 & \Delta_2 &          &          \\
                     & \ddots   & \ddots   &          \\
                     &          & \Gamma_k & \Delta_k \\
            \end{bmatrix}\!,
    \mathhfill
    \Delta_k = \begin{bmatrix}
                 1        &   \\
                 \delta_k & 1 \\
               \end{bmatrix}\!,
    \mathhfill
    \Gamma_k = \begin{bmatrix}
                 & \sigma_k  \\
                 \eta_k & \lambda_k
               \end{bmatrix}\!.
  \]
  If we initialize $d_{-1} = d_0 = \sigma_1 = \eta_1 = \lambda_1 = 0$, individual factorization steps are obtained from the recursion formulae
  \begin{subequations}
    \label{ldlt_formula}
      \begin{alignat}{2}
        d_{2k-1} &= 1 - \sigma_k^2 d_{2k-2}, \quad & k & \ge 1, \\
        d_{2k} &= -1 - \eta_k^2 d_{2k-3} - \lambda_k^2 d_{2k-2} - \delta_k^2 d_{2k-1}, \quad & k & \ge 1, \\
        \delta_k &= (\alpha_k - \lambda_k \sigma_k d_{2k-2}) / d_{2k-1}, \quad & k & \ge 1, \\
        \sigma_k &=  \beta_k / d_{2k-2}, \quad & k & \ge 2, \\
        \eta_k &= \gamma_k / d_{2k-3}, \quad & k & \ge 2, \\
        \lambda_k &= -\eta_k \delta_{k-1} d_{2k-3} / d_{2k-2}, \quad & k & \ge 2.
    \end{alignat}
  \end{subequations}

  \subsubsection{Update of the \TriCG iterate}

  In order to compute the solution $z_k$ of~\eqref{eq:sub-tricg}, we update the solution $p_k := (\pi_1, \cdots, \pi_{2k})$ of $L_k D_k p_k = (\beta_1 e_1 + \gamma_1 e_2)$.
  The components of $p_k$ are computed from

  \begin{subequations}
    \label{pk}
    \begin{align}
      \pi_{2k-1} &=
      \begin{cases}
        \beta_1 / d_1, & k = 1, \\
        -\sigma_k d_{2k-2} \pi_{2k-2} / d_{2k-1}, & k \ge 2,
      \end{cases}
      \\
      \pi_{2k} &=
      \begin{cases}
        (\gamma_1 - \delta_1 \beta_1) / d_2, & k = 1, \\
        -(\delta_k d_{2k-1} \pi_{2k-1} + \lambda_k d_{2k-2} \pi_{2k-2} + \eta_k d_{2k-3} \pi_{2k-3}) / d_{2k}, & k \ge 2.
      \end{cases}
    \end{align}
  \end{subequations}
  If we were to update $(x_k,y_k)$ directly from~\eqref{eq:xy}, all components of $z_k := (\zeta_1, \cdots, \zeta_{2k})$ would have to be recomputed because of the backward substitution required to solve $L_k^T z_k = p_k$, which would require us to store $W_k$ entirely.
  To avoid such drawbacks, we employ the strategy of \cite{paige-saunders-1975}.
  Let
  \begin{equation}
    \label{eq:def-Gk}
  G_k := W_k L_k^{-T} \; \Longleftrightarrow \; L_k G_k^T = W_k^T,~G_k = \begin{bmatrix} G^x_k \\ G^y_k \end{bmatrix} = \begin{bmatrix} g^x_1 & \cdots & g^x_{2k} \\ g^y_1 & \cdots & g^y_{2k} \end{bmatrix},
  \end{equation}
  defined by $g^x_{-1} = g^x_0 = g^y_{-1} = g^y_0 = 0$, and the recursion
  \begin{equation}
    \begin{aligned}
      g^x_{2k-1} &= - \sigma_k g^x_{2k-2} + v_k , \\
      g^y_{2k-1} &= - \sigma_k g^y_{2k-2}, \\
      g^x_{2k}   &= - \delta_k g^x_{2k-1} - \lambda_k g^x_{2k-2} - \eta_k g^x_{2k-3}, \\
      g^y_{2k}   &= - \delta_k g^y_{2k-1} - \lambda_k g^y_{2k-2} - \eta_k g^y_{2k-3} + u_k.
    \end{aligned}
  \end{equation}
  This gives $(x_k,y_k) = W_k z_k = G_k L_k^T z_k = G_k p_k$ and the solution may be updated efficiently as
  \begin{subequations}
    \label{eq:tricg}
    \begin{align}
      x_k &= G^x_k p_k = x_{k-1} + \pi_{2k-1} g^x_{2k-1} + \pi_{2k} g^x_{2k}, \\
      y_k &= G^y_k p_k = y_{k-1} + \pi_{2k-1} g^y_{2k-1} + \pi_{2k} g^y_{2k}.
    \end{align}
  \end{subequations}

  \subsubsection{Residual computation}

  The expression~\eqref{eq:residual_general} combines with~\eqref{eq:sub-tricg} to yield the residual at the \TriCG iterate:
  \begin{align}
    r_k & =
    -H W_k (S_k z_k - \beta_1 e_1 - \gamma_1 e_2 ) - H w_{k+1} \Psi_{k+1}^T \bmat{e_{2k-1}^T \\ e_{2k}^T} z_k
    \nonumber \\ & =
    - H w_{k+1} \Psi_{k+1}^T \bmat{\zeta_{2k-1} \\ \zeta_{2k}}
    \nonumber \\ & =
    - H w_{k+1} \bmat{\beta_{k+1} \zeta_{2k} \\ \gamma_{k+1} \zeta_{2k-1}}.
    \label{eq:residual_tricg}
  \end{align}
  Because $L_k^T z_k = p_k$, we have $\zeta_{2k} = \pi_{2k}$ and $\zeta_{2k-1} = \pi_{2k-1} - \delta_k \pi_{2k}$.
  Therefore, with the relation $w_{k+1}^T H w_{k+1} = I_2$, it is natural to measure the residual in the \(H^{-1}\)-norm:
  \begin{subequations}
    \label{eq:||r||_tricg}
    \begin{align}
      \|r_0\|_{H^{-1}} &= \sqrt{\beta_1^2 + \gamma_1^2},
      \label{eq:||r0||_tricg} \\
      \|r_k\|_{H^{-1}} &= \sqrt{\gamma_{k+1}^2 (\pi_{2k-1} - \delta_k \pi_{2k})^2 + \beta_{k+1}^2 \pi_{2k}^2}, \quad k \ge 1.
      \label{eq:||rk||_tricg}
    \end{align}
  \end{subequations}
  We summarize the complete procedure as \Cref{alg:tricg}.

  \begin{algorithm}[t]
    \caption{%
      \TriCG
    }
    \label{alg:tricg}
    \begin{algorithmic}[1]
      \Require $A$, $b$, $c$, $v \mapsto M \backslash v$, $u \mapsto N \backslash u$
      \State $x_0 = 0$, $y_0 = 0$
      \State $g^x_{-1} = 0$, $g^x_0 = 0$, $g^y_{-1} = 0$, $g^y_0 = 0$
      \State $u_0 = 0$, $v_0 = 0$ \Comment{begin orthogonal triorthogonalization}
      \State $\beta_1 M v_1 = b$, $\gamma_1 N u_1 = c$ \Comment{$(\beta_1, \gamma_1) > 0$ so that $\|v_1\|_M = \|u_1\|_N = 1$}
      \State $\|r_0\|_{H^{-1}} = (\beta_1^2 + \gamma_1^2)^{\frac12}$ \Comment{compute $\|r_0\|_{H^{-1}}$}
      \State $d_{-1} = d_0 = \sigma_1 = \eta_1 = \lambda_1 = 0$ \Comment{initialize the $LDL^T$ factorization}
      \For{\(k = 1, 2, \dots\)}
        \State $q = A u_k - \gamma_k M v_{k-1}$, $\alpha_k = v_k^T q$ \Comment{continue orthogonal triorthogonalization}
        \State $p = A^T v_k - \beta_k N u_{k-1}$
        \State $\beta_{k+1} M v_{k+1} = q - \alpha_k M v_k$ \Comment{$\beta_{k+1} > 0$ so that $\|v_{k+1}\|_M = 1$}
        \State $\gamma_{k+1} N u_{k+1} = p - \alpha_k N u_k$ \Comment{$\gamma_{k+1} > 0$ so that $\|u_{k+1}\|_N = 1$}
        \State $d_{2k-1} = 1 - \sigma_k^2 d_{2k-2}$ \Comment{continue the $LDL^T$ factorization}
        \State $\delta_k = (\alpha_k - \lambda_k \sigma_k d_{2k-2}) / d_{2k-1}$ \Comment{compute $\Delta_k$}
        \State $d_{2k} = -1 - \eta_k^2 d_{2k-3} - \lambda_k^2 d_{2k-2} - \delta_k^2 d_{2k-1}$ \Comment{update $D_k$}
        \If{$k == 1$}
          \State $\pi_{2k-1} = \beta_k / d_{2k-1}$
          \State $\pi_{2k} = (\gamma_k - \delta_k \beta_k) / d_{2k}$ \Comment{initial solution of $L_k D_k p_k = \beta_1 e_1 + \gamma_1 e_2$}
        \Else
          \State $\sigma_k = \beta_k / d_{2k-2}$
          \State $\eta_k = \gamma_k / d_{2k-3}$ \Comment{compute $\Gamma_k$}
          \State $\lambda_k = - (\eta_k \delta_{k-1} d_{2k-3}) / d_{2k-2}$
          \State $\pi_{2k-1} = -(\sigma_k \pi_{2k-2} d_{2k-2}) / d_{2k-1}$ \Comment{update $p_k$}
          \State $\pi_{2k} = -(\delta_k \pi_{2k-1} d_{2k-1} + \lambda_k \pi_{2k-2} d_{2k-2} + \eta_k \pi_{2k-3} d_{2k-3}) / d_{2k}$
        \EndIf
        \State $g^x_{2k-1} = v_k - \sigma_k g^x_{2k-2}$ \Comment{update $G_k^x$}
        \State $g^x_{2k} = - \delta_k g^x_{2k-1} - \lambda_k g^x_{2k-2} - \eta_k g^x_{2k-3}$
        \State $g^y_{2k-1} = - \sigma_k g^y_{2k-2}$ \Comment{update $G_k^y$}
        \State $g^y_{2k} = u_k - \delta_k g^y_{2k-1} - \lambda_k g^y_{2k-2} - \eta_k g^y_{2k-3}$
        \State $x_k = x_{k-1} + \pi_{2k-1} g^x_{2k-1} + \pi_{2k} g^x_{2k}$ \Comment{update $x_k$}
        \State $y_k = y_{k-1} + \pi_{2k-1} g^y_{2k-1} + \pi_{2k} g^y_{2k}$ \Comment{update $y_k$}
        \State $\|r_k\|_{H^{-1}} = (\gamma_{k+1}^2 (\pi_{2k-1} - \delta_k \pi_{2k})^2 + \beta_{k+1}^2 \pi_{2k}^2)^{\frac12}$ \Comment{compute $\|r_k\|_{H^{-1}}$}
      \EndFor
    \end{algorithmic}
  \end{algorithm}

  \subsubsection{Storage}

  \TriCG requires one operator-vector product with $A$ and one with $A\T$ per iteration.
  With the assumption that in-place \emph{gemv} updates of the form $y \leftarrow Au + \gamma y$ and $y \leftarrow A\T v + \beta y$ are available, \TriCG requires five \(n\)-vectors ($y_k$, $u_{k-1}$, $u_k$, $g^y_{2k-1}$, $g^y_{2k}$) and five \(m\)-vectors ($x_k$, $v_{k-1}$, $v_k$, $g^x_{2k-1}$, $g^x_{2k}$).
  If in-place \emph{gemv} updates are not available, additional $m$- and $n$-vectors are required to store \(Au\) and \(A\T v\).
  Note that $A$, $A\T$, $M^{-1}$ and $N^{-1}$ do not need to be formed explicitly, and can be implemented as abstract operators.
  For instance, we could compute the Cholesky factorization of $M$ and $N$ and create abstract operators that perform the forward and backsolves.
  Extra $m$- and $n$-vectors could be necessary to store the results of those operators.

  \subsubsection{Properties}

  In this section, we formulate optimality properties of the \TriCG iterates.

  \begin{proposition}
    \label{prop:min-max}
    The \(k\)-th \TriCG iterate $(x_k,y_k)$ solves
    \begin{equation}
    \label{eq:min-max}
     \minimize{x \in \R^m} \, \maximize{y \in \R^n} \, \mathcal{L}(x,y)
      \quad \st \, \begin{bmatrix} x \\ y \end{bmatrix} \in \Range(W_k),
    \end{equation}
    where $\mathcal{L}(x,y) = \tfrac{1}{2}\|x\|_M^2 - \tfrac{1}{2}\|y\|_N^2 + x\T A y - b\T x - c\T y$.
    Equivalently, $(x_k,y_k)$ solves
    \begin{equation}
      \label{eq:error-min-max}
      \minimize{x \in \R^m} \, \maximize{y \in \R^n} \, \mathcal{E}(x,y)
      \quad \st \, \begin{bmatrix} x \\ y \end{bmatrix} \in \Range(W_k),
    \end{equation}
    where \(\mathcal{E}(x,y)\) is the indefinite error metric
    \[
      \mathcal{E}(x,y) :=
      e_{r}^T
      \begin{bmatrix}
        M & \phantom{-}A \\
        A^T & -N
      \end{bmatrix}
      e_{r},
      \qquad
      e_{r} := (x^{\ast} - x,~y^{\ast} - y),
    \]
    and $(x^{\ast},y^{\ast})$ is the exact solution of~\eqref{eq:sqd}.
  \end{proposition}

  \begin{proof}
    $\mathcal{L}(x,y)$ is strictly convex in $x$ because $\nabla_{xx}^2 \mathcal{L}(x,y) = M \succ 0$ and strictly concave in $y$ because $\nabla_{yy}^2 \mathcal{L}(x,y) = -N \prec 0$.
    Therefore,~\eqref{eq:min-max} admits a unique solution because the feasible set $\Range(W_k) \neq  \varnothing$.
    Its first-order optimality conditions are
    \[
      W_k^T
      \begin{bmatrix}
        Mx + Ay - b \\
        A\T x -Ny -c
      \end{bmatrix}
      = 0,
    \]
    and coincide with~\eqref{eq:galerkin-tricg}.
    The rest of the proof follows from the fact that \(\mathcal{L}(x,y)\) and \(\mathcal{E}(x,y)\) are equal up to a constant.
  \end{proof}

\begin{figure}[ht]
    \centering
    \includetikzgraphics[width=0.7\textwidth]{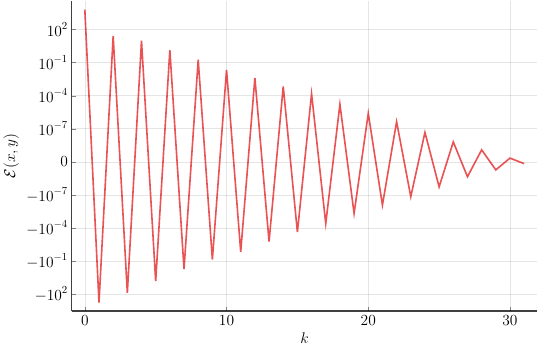}
    \label{fig:metric}
    \caption{$\mathcal{E}(x,y)$ history of \TriCG on illc1850, where $A$ is $1850 \times 712$ with $8636$ nonzeros.}
  \end{figure}

  Although SQD matrices are indefinite, $\mathcal{E}(x,y)$ can be seen as a metric that generalizes the energy norm.
  A similar metric is used by \cite{orban-arioli-2017} in the context of their \emph{generalized conjugate gradient method} for SQD systems.
  \Cref{fig:metric} illustrates the evolution of \(\mathcal{E}(x_k, y_k)\) along the \TriCG iterations on problem \emph{illc1850}, to be described in \cref{sec:num-results}, where oscillations from positive to negative values and decreasing amplitude are evident.

  \subsection{Derivation of \TriMR}

  In the same way as \TriCG is related to block-\CGM, the minimum residual variant \TriMR developed below is related to block-\MINRES.
  The $k$-th \TriMR iterate is defined as the solution of the linear least-squares problem
  \begin{equation}
    \label{eq:sub-trimr}
    \minimize{(x_k,y_k) \in \Range(W_k)} \|r_k\|_{H^{-1}}
    \quad \Longleftrightarrow \quad
    \minimize{z_k \in \R^{2k}} \|S_{k+1,k} z_k - \beta_1 e_1 - \gamma_1 e_2\|,
  \end{equation}
  where the equivalence follows from~\eqref{eq:residual_general}.
  We now outline the main stages of the subproblem solution.

  \subsubsection{A QR factorization}

  The solution of~\eqref{eq:sub-trimr} can be determined via the QR factorization
  \begin{equation}
    \label{eq:qr_Sk}
    S_{k+1, k} = Q_k \begin{bmatrix} R_k \\ 0 \end{bmatrix},
  \end{equation}
  which can be updated at each iteration, where $Q_k \in \R^{(2k+2) \times (2k+2)}$ is a product of Givens reflections, and
  \begin{equation}
    \label{eq:Rk}
    R_k = \begin{bmatrix}
             \delta_1 & \sigma_1 & \eta_1 & \lambda_1 & \mu_1  &        & \\
                      & \ddots   & \ddots & \ddots    & \ddots & \ddots & \\
                      &          & \ddots & \ddots    & \ddots & \ddots & \mu_{2k-4} \\
                      &          &        & \ddots    & \ddots & \ddots & \lambda_{2k-3} \\
                      &          &        &           & \ddots & \ddots & \eta_{2k-2} \\
                      &          &        &           &        & \ddots & \sigma_{2k-1} \\
                      &          &        &           &        &        & \delta_{2k}
          \end{bmatrix} \in \R^{(2k) \times (2k)}.
  \end{equation}
  Below, we outline the main steps of the update and commit all details to \cref{sec:appendix} for reference.

  At iteration \(k\), four reflections are necessary to update~\eqref{eq:qr_Sk}.
  Together, they affect four rows and six columns of \(S_{k+1,k}\).
  We denote their product $Q_{2k-1,2k+2}$ in \eqref{eq:Q} so that $Q_k^T = Q_{2k-1,2k+2} \cdots Q_{1,4}$.
  If we initialize $\bar{\theta}_1 := \alpha_1$, $\bar{\delta}_1 := 1$, $\bar{\delta}_2 := -1$, $\bar{\sigma}_1 := \alpha_1$, $\bar{\eta}_1 := 0$, $\bar{\lambda}_1 := \gamma_2$ and $\bar{\sigma}_2 := \beta_2$, individual factorization steps may be represented as an application of $Q_{2k-1,2k+2}$ to $Q_{k-1}^T S_{k+1,k}$:
  \[
    \kbordermatrix{
         & 2k-1                & 2k                  &        & 2k+1               & 2k+2                 & 2k+3        & 2k+4 \\
    2k-1 & \bar{\delta}_{2k-1} & \bar{\sigma}_{2k-1} & \vrule & \bar{\eta}_{2k-1}  & \bar{\lambda}_{2k-1} & 0           & 0 \\
    2k   & \bar{\theta}_{k}    & \bar{\delta}_{2k}   & \vrule &  \bar{\sigma}_{2k} & 0                    & 0           & 0 \\
    2k+1 & 0                   & \gamma_{k+1}        & \vrule &  1                 & \alpha_{k+1}         & 0           & \gamma_{k+2} \\
    2k+2 & \beta_{k+1}         & 0                   & \vrule &  \alpha_{k+1}      & -1                   & \beta_{k+2} & 0
    }.
  \]
  Because $\alpha_{k+1}$, $\gamma_{k+2}$ and $\beta_{k+2}$ are not yet available at iteration $k$, we apply the last four reflections at iteration $k+1$ to determine all remaining coefficients of rows $2k-1$ and $2k$ of $R_{k+1}$ and $R_{k+2}$:
  \[
    \kbordermatrix{
         & 2k-1          & 2k            &        & 2k+1                & 2k+2                & 2k+3                & 2k+4 \\
    2k-1 & \delta_{2k-1} & \sigma_{2k-1} & \vrule & \eta_{2k-1}         & \lambda_{2k-1}      & \mu_{2k-1}          & 0                    \\
    2k   & 0             & \delta_{2k}   & \vrule & \sigma_{2k}         & \eta_{2k}           & \lambda_{2k}        & \mu_{2k}             \\
    2k+1 & 0             & 0             & \vrule & \bar{\delta}_{2k+1} & \bar{\sigma}_{2k+1} & \bar{\eta}_{2k+1}   & \bar{\lambda}_{2k+1} \\
    2k+2 & 0             & 0             & \vrule & \bar{\theta}_{k+1}  & \bar{\delta}_{2k+2} & \bar{\sigma}_{2k+2} & 0
    }.
  \]
  Additional details about the four reflections that compose $Q_{2k-1,2k+4}$ and the factorization are available in \eqref{eq:ref1}--\eqref{eq:ref4}.

  \subsection{Update of the \TriMR iterate}

  We can avoid computing $z_k$ in~\eqref{eq:sub-trimr} for the same reasons as with \TriCG by updating instead $p_k := (\pi_1, \cdots, \pi_{2k})$, which is defined via $R_k z_k = p_k$:
  \begin{subequations}
    \label{eq:pk}
    \begin{align}
      \bar{p}_0 &:= (\beta_1, \gamma_1), \\
      \bar{p}_k &:=
      (p_k, \bar{\pi}_{2k+1}, \bar{\pi}_{2k+2}) =
      Q_k^T (\beta_1 e_1 + \gamma_1 e_2), \quad k \ge 1.
    \end{align}
  \end{subequations}
  $\bar{p}_{k}$ can be easily determined from $\bar{p}_{k-1}$ because $\bar{p}_{k} = Q_{2k-1,2k+2}~(\bar{p}_{k-1}, 0, 0)$.
  Details are given in \eqref{eq:ref12} and \eqref{eq:ref34}.
  We set
  \begin{equation}
    \label{eq:Gk-trimr}
    G_k := W_k R_k^{-1} \quad \Longleftrightarrow \quad R_k G_k^T = W_k^T,
  \end{equation}
  similarly to~\eqref{eq:def-Gk}.
  The columns of \(G_k\) are obtained from the recursion
  \begin{equation}
    \label{eq:def-Gk2}
    \begin{aligned}
      g^x_{2k-1} &= (v_k           - \mu_{2k-5} g^x_{2k-5} - \lambda_{2k-4} g^x_{2k-4} - \eta_{2k-3} g^x_{2k-3} - \sigma_{2k-2} g^x_{2k-2}) / \delta_{2k-1} \\
      g^y_{2k-1} &= (\phantom{u_k} - \mu_{2k-5} g^y_{2k-5} - \lambda_{2k-4} g^y_{2k-4} - \eta_{2k-3} g^y_{2k-3} - \sigma_{2k-2} g^y_{2k-2}) / \delta_{2k-1} \\
      g^x_{2k}   &= (\phantom{v_k} - \mu_{2k-4} g^x_{2k-4} - \lambda_{2k-3} g^x_{2k-3} - \eta_{2k-2} g^x_{2k-2} - \sigma_{2k-1} g^x_{2k-1}) / \delta_{2k}   \\
      g^y_{2k}   &= (u_k           - \mu_{2k-4} g^y_{2k-4} - \lambda_{2k-3} g^y_{2k-3} - \eta_{2k-2} g^y_{2k-2} - \sigma_{2k-1} g^y_{2k-1}) / \delta_{2k},
    \end{aligned}
  \end{equation}
  where we set $\eta_j$, $\lambda_j$, $\mu_j$, $g^x_j$ and $g^y_j$ to zero if $j \le 0$.
  Analogously to \TriCG, $(x_k,y_k) = W_k z_k = G_k R_k z_k = G_k p_k$ and the solution may be updated efficiently as \eqref{eq:tricg}.

  \subsubsection{Residual computation}

  The definition of \(\bar{p}_k\),~\eqref{eq:residual_general} and~\eqref{eq:qr_Sk} yield
  \begin{equation}
    \label{eq:residual_trimr}
    \|r_k\|_{H^{-1}} =
    \|S_{k+1,k} z_k - (\beta_1 e_1 + \gamma_1 e_2)\| =
    \left\| \bmat{R_k \\ 0} z_k - \bar{p}_k \right\| =
    \sqrt{\bar{\pi}_{2k+1}^2 + \bar{\pi}_{2k+2}^2}.
  \end{equation}

  The complete algorithm is stated as \Cref{alg:trimr}.

  \begin{algorithm}[t]
    \caption{%
      \TriMR
    }
    \label{alg:trimr}
    \begin{algorithmic}[1]
      \Require $A$, $b$, $c$, $v \mapsto M \backslash v$, $u \mapsto N \backslash u$
      \State $x_0 = 0$, $y_0 = 0$
      \State $g^x_{-3} = 0$, $g^x_{-2} = 0$, $g^x_{-1} = 0$, $g^x_0 = 0$
      \State $g^y_{-3} = 0$, $g^y_{-2} = 0$, $g^y_{-1} = 0$, $g^y_0 = 0$
      \State $u_0 = 0$, $v_0 = 0$ \Comment{begin orthogonal triorthogonalization}
      \State $\beta_1 M v_1 = b$, $\gamma_1 N u_1 = c$ \Comment{$(\beta_1, \gamma_1) > 0$ so that $\|v_1\|_M = \|u_1\|_N = 1$}
      \State $\|r_0\|_{H^{-1}} = (\beta_1^2 + \gamma_1^2)^{\frac12}$ \Comment{compute $\|r_0\|_{H^{-1}}$}
      %\State $\bar{\pi}_1 = \beta_1$, $\bar{\pi}_2 = \gamma_1$ \Comment{Initialize $\bar{p}_0$}
      \For{\(k = 1, 2, \dots\)}
        \State $q = A u_k - \gamma_k M v_{k-1}$, $\alpha_k = v_k^T q$ \Comment{continue orthogonal triorthogonalization}
        \State $p = A^T v_k - \beta_k N u_{k-1}$
        \State $\beta_{k+1} M v_{k+1} = q - \alpha_k M v_k$ \Comment{$\beta_{k+1} > 0$ so that $\|v_{k+1}\|_M = 1$}
        \State $\gamma_{k+1} N u_{k+1} = p - \alpha_k N u_k$ \Comment{$\gamma_{k+1} > 0$ so that $\|u_{k+1}\|_N = 1$}
        \If{$k == 1$}
          \State $\bar{\theta}_1 = \alpha_1$, $\bar{\delta}_1 = 1$, $\bar{\delta}_2 = -1$ \Comment{initialize the $QR$ factorization}
          \State $\bar{\sigma}_1 = \alpha_1$, $\bar{\eta}_1 = 0$, $\bar{\lambda}_1 = \beta_2$, $\bar{\sigma}_2 = \gamma_2$
        \Else
          \State Compute $\eta_{2k-3}$, $\lambda_{2k-3}$, $\mu_{2k-3}$, $\sigma_{2k-2}$, $\eta_{2k-2}$, $\lambda_{2k-2}$, $\mu_{2k-2}$ \Comment{update $R_k$}
          \State Compute $\bar{\theta}_k$, $\bar{\delta}_{2k-1}$, $\bar{\delta}_{2k}$, $\bar{\sigma}_{2k-1}$, $\bar{\eta}_{2k-1}$, $\bar{\lambda}_{2k-1}$, $\bar{\sigma}_{2k}$
        \EndIf
        \State Compute $Q_{2k-1,2k+2}$, $\delta_{2k-1}$, $\sigma_{2k-1}$, $\delta_{2k}$ \Comment{continue the $QR$ factorization}
        \State Compute $\pi_{2k-1}$, $\pi_{2k}$, $\bar{\pi}_{2k+1}$, $\bar{\pi}_{2k+2}$ \Comment{update $\bar{p}_k$}
        \State {\small $g^x_{2k-1} = (v_k - \mu_{2k-5} g^x_{2k-5} - \lambda_{2k-4} g^x_{2k-4} - \eta_{2k-3} g^x_{2k-3} - \sigma_{2k-2} g^x_{2k-2}) / \delta_{2k-1}$}
        \State {\small $g^x_{2k} = -(\mu_{2k-4} g^x_{2k-4} + \lambda_{2k-3} g^x_{2k-3} + \eta_{2k-2} g^x_{2k-2} + \sigma_{2k-1} g^x_{2k-1}) / \delta_{2k}$} \Comment{update $G_k^x$}
        \State {\small $g^y_{2k-1} = -(\mu_{2k-5} g^y_{2k-5} + \lambda_{2k-4} g^y_{2k-4} + \eta_{2k-3} g^y_{2k-3} + \sigma_{2k-2} g^y_{2k-2}) / \delta_{2k-1}$}
        \State {\small $g^y_{2k} = (u_k - \mu_{2k-4} g^y_{2k-4} - \lambda_{2k-3} g^y_{2k-3} - \eta_{2k-2} g^y_{2k-2} - \sigma_{2k-1} g^y_{2k-1}) / \delta_{2k}$} \Comment{update $G_k^y$}
        \State $x_k = x_{k-1} + \pi_{2k-1} g^x_{2k-1} + \pi_{2k} g^x_{2k}$ \Comment{update $x_k$}
        \State $y_k = y_{k-1} + \pi_{2k-1} g^y_{2k-1} + \pi_{2k} g^y_{2k}$ \Comment{update $y_k$}
        \State $\|r_k\|_{H^{-1}} = (\bar{\pi}_{2k+1}^2 + \bar{\pi}_{2k+2}^2)^{\frac12}$ \Comment{compute $\|r_k\|_{H^{-1}}$}
      \EndFor
    \end{algorithmic}
  \end{algorithm}

  \subsubsection{Storage}

  \TriMR has the same storage requirements as \TriCG plus two $n$-vectors ($g^y_{2k-2}$, $g^y_{2k-3}$) and two $m$-vectors ($g^x_{2k-2}$, $g^x_{2k-3}$).
  All other vectors are identical to those in \TriCG.

  \section{Implementation and numerical experiments}
  \label{sec:num-results}

  We evaluate the performance of \TriCG and \TriMR on SQD systems generated from rectangular matrices \(A\) obtained from the UFL collection of \cite{davis-hu-2011}.\footnote{Now the SuiteSparse Matrix Collection \https{sparse.tamu.edu}.}
  We implemented \Cref{alg:tricg} and \Cref{alg:trimr} in Julia\footnote{\https{julialang.org}} \citep{bezanson-edelman-karpinski-shah-2017}, version \(1.5\).
  Both algorithms are available as part of the \texttt{Krylov.jl} collection of Krylov methods \citep{montoison-orban-krylov-2020}.

  Because standard \CGM may break down when applied to~\eqref{eq:sqd}, we compare the evolution of the \TriCG residual to that of \SYMMLQ, whose iterates are always well defined.
  % If the \CGM iterate is well defined at a given \SYMMLQ iteration, we step to it with the help of an inexpensive transfer procedure described by \cite{paige-saunders-1975}.
  Similarly, we compare the evolution of the \TriMR residual to that of \MINRES.
  In order to evaluate benefits of \TriCG and \TriMR in terms of loss of orthogonality along the iterations, we also compare the evolution of \TriCG and \TriMR residuals to those of \BLCG and \BLMINRES, respectively, applied to $K$ with block right-hand side $B$ where the two approximate solutions are summed at the last iteration.
  \SYMMLQ, \BLCG, \MINRES and \BLMINRES are run with preconditioner \(H\).

  In our first set of experiments, we set \(M\) and \(N\) to the identity.
  Thus the \(H^{-1}\)-norm is simply the Euclidean norm.
  The right-hand side \((b, c)\) is generated such that the exact solution of~\eqref{eq:sqd} is the vector of ones.
  Residuals $r_k = b - Ax_k$ are calculated explicitly at each iteration in order to evaluate $\|r_k\|$ instead of using~\eqref{eq:residual_tricg} or~\eqref{eq:residual_trimr}.
  Each algorithm stops as soon as $\|r_k\| \leq \varepsilon_a + \|(b,c)\| \varepsilon_r$ with absolute tolerance $\varepsilon_a = 10^{-12}$ and relative tolerance $\varepsilon_r = 10^{-10}$.

\begin{figure}[ht]
    \centering
    \includetikzgraphics[width=0.49\textwidth]{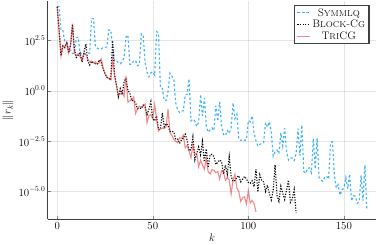}
    \hfill
    \includetikzgraphics[width=0.49\textwidth]{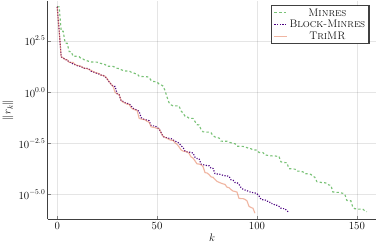}
    \label{fig:lp_czprob}
    \caption{Residual history of \SYMMLQ, \BLCG, \TriCG, \MINRES, \BLMINRES and \TriMR on lp\_czprob, where $A$ is $929 \times 3562$ with $10708$ nonzeros.}
  \end{figure}

\begin{figure}[ht]
    \centering
    \includetikzgraphics[width=0.49\textwidth]{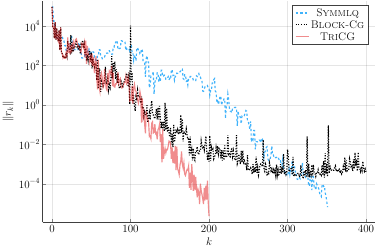}
    \hfill
    \includetikzgraphics[width=0.49\textwidth]{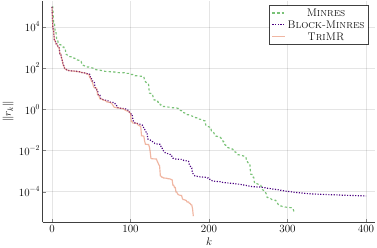}
    \label{fig:lp_osa_07}
    \caption{Residual history of \SYMMLQ, \BLCG, \TriCG, \MINRES, \BLMINRES and \TriMR on lp\_osa\_07, where $A$ is $1108 \times 25067$ with $144812$ nonzeros.}
  \end{figure}

  \Cref{fig:lp_osa_07,fig:lp_czprob} report residual histories on matrices arising from linear optimization.
  In all cases, the \TriCG and \TriMR residuals attain the required tolerance in around half the number of iterations of \SYMMLQ and \MINRES, respectively.
  We also note that the \TriCG and \BLCG residuals are close, but not quite superperposed, as are the \MINRES and \BLMINRES residuals.
  These results are encouraging if ones wishes to employ \TriCG or \TriMR to solve the linear systems arising at each iteration of a numerical method for constrained optimization, including interior-point methods, where the systems have the form of those just tested.
  The results also suggest that orthogonality is not lost quite as fast in \TriCG and \TriMR as it is in \BLCG and \BLMINRES.

\begin{figure}[ht]
    \centering
    \includetikzgraphics[width=0.49\textwidth]{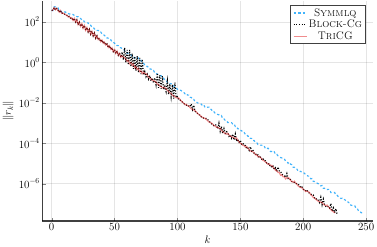}
    \hfill
    \includetikzgraphics[width=0.49\textwidth]{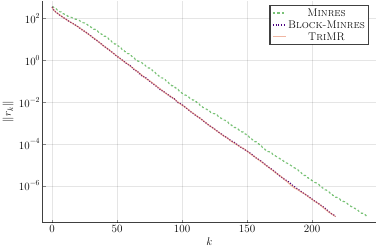}
    \label{fig:Maragal_6}
    \caption{Residual history of \SYMMLQ, \BLCG, \TriCG, \MINRES, \BLMINRES and \TriMR on Maragal\_6, where $A$ is $21255 \times 10152$ with $537694$ nonzeros.}
  \end{figure}

\begin{figure}[ht]
    \centering
    \includetikzgraphics[width=0.49\textwidth]{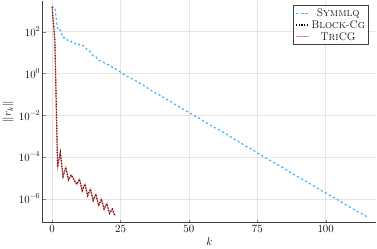}
    \hfill
    \includetikzgraphics[width=0.49\textwidth]{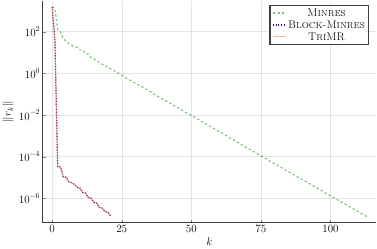}
    \label{fig:landmark}
    \caption{Residual history of \SYMMLQ, \BLCG, \TriCG, \MINRES, \BLMINRES and \TriMR on landmark, where $A$ is $71952 \times 2704$ with $1146848$ nonzeros.}
  \end{figure}

  \Cref{fig:Maragal_6,fig:landmark} report residual histories on matrices arising from least-squares problems.
  In all cases, \TriCG and \TriMR require fewer iterations than \SYMMLQ and \MINRES.
  On these two problems, the residuals of the block methods are nearly superposed.
  We observe on our test problems that \TriCG and \TriMR perform fewer iterations when the singular values of $A$ are clustered.
  However, a deeper analysis is required to confirm this empirical observation.

  We simulate the behavior of the six methods in exact arithmetic in hopes to compare loss of orthogonality empirically.
  \Cref{fig:lp_osa_07_mp} reports residual histories obtained when running each method on lp\_osa\_07 entirely in quadruple and octuple precision, and should be compared with \Cref{fig:lp_osa_07}.
  This time, the \TriCG and \BLCG residuals are nearly superposed, as are the \TriMR and \BLMINRES residuals as the increased accuracy of arithmetic mitigates loss of orthogonality.
  The block methods continue to require around half as many iterations as the standard methods.
  However, we note an additional phenomenon: doubling the number of digits reduces the number of iterations by a factor of approximately two.

  \begin{figure}[ht]
    \centering
    \includetikzgraphics[width=0.49\textwidth]{mp/lp_osa_07_tricg_blcg_symmlq_Float128}
    \hfill
    \includetikzgraphics[width=0.49\textwidth]{mp/lp_osa_07_trimr_blminres_minres_Float128}
    \vfill
    \includetikzgraphics[width=0.49\textwidth]{mp/lp_osa_07_tricg_blcg_symmlq_Float256}
    \hfill
    \includetikzgraphics[width=0.49\textwidth]{mp/lp_osa_07_trimr_blminres_minres_Float256}
    \label{fig:lp_osa_07_mp}
    \caption{Residual history of \SYMMLQ, \BLCG, \TriCG, \MINRES, \BLMINRES and \TriMR iterates on lp\_osa\_07 in quadruple (top) and octuple precision (bottom).}
  \end{figure}

  In a second set of experiments, we run all six methods on discretized and stabilized Stokes equations generated by the MATLAB package
  % \emph{Incompressible Flow Iterative Solution Software} (IFISS)
  IFISS, version~\(3.6\), of \cite{ifiss}.
  Whenever the discrete velocity and pressure belong to finite-element spaces that do not satisfy the \textit{inf-sup}, or Ladyzhenskaya-Babuška-Brezzi (LBB), stability conditions \citep{boffi-brezzi-fortin-2013}, a nonzero and negative semi-definite stabilization term $-N$ is inserted in the bottom block of~\eqref{eq:saddle_point}.
  It is the case with the unstable 2D finite-element pairs $Q_1$-$P_0$ and $Q_1$-$Q_1$, which we use on a test problem from IFISS.
  In order to obtain an SQD system, we add $10^{-5}I$ to $N$.
  For this set of problems, $M$ and $N$ are not identity operators, and each algorithm stops as soon as $\|r_k\|_{H^{-1}} \leq \varepsilon_a + \|(b,c)\|_{H^{-1}} \varepsilon_r$ with the same tolerances as above.

  \begin{figure}[ht]
    \centering
    \includetikzgraphics[width=0.49\textwidth]{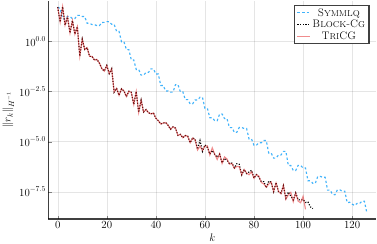}
    \hfill
    \includetikzgraphics[width=0.49\textwidth]{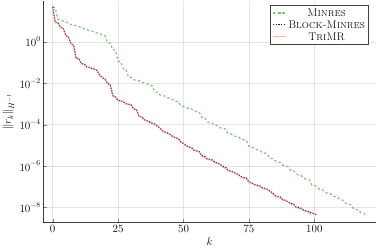}
    \label{fig:channel_domain1}
    \caption{Residual history of \SYMMLQ, \BLCG, \TriCG, \MINRES, \BLMINRES and \TriMR on channel\_domain problem with \(Q_1\)-\(P_0\) discretization.
    The discretized linear system has size $12546 \times 12546$ with $147742$ nonzeros.}
  \end{figure}

  \begin{figure}[ht]
    \centering
    \includetikzgraphics[width=0.49\textwidth]{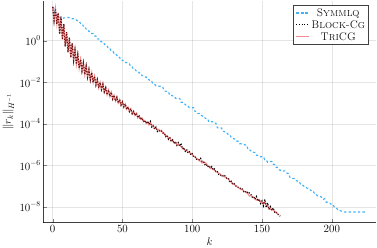}
    \hfill
    \includetikzgraphics[width=0.49\textwidth]{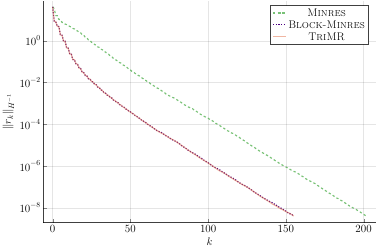}
    \label{fig:channel_domain2}
    \caption{Residual history of \SYMMLQ, \BLCG, \TriCG, \MINRES, \BLMINRES and \TriMR on channel\_domain problem with \(Q_1\)-\(Q_1\) discretization.
    The discretized linear system has size $12675 \times 12675$ with $242381$ nonzeros.}
  \end{figure}

  \Cref{fig:channel_domain1,fig:channel_domain2} report residual histories.
  \TriCG and \TriMR reach the prescribed tolerance before \SYMMLQ and \MINRES with around 25\% fewer iteration.
  These results suggest that \TriCG and \TriMR may be of interest to solve other PDEs whose discretization leads to~\eqref{eq:sqd}, such as the Reissner-Mindlin plate model in linear elasticity \citep{braess-2007}.

  \subsection{Discussion and extensions}

  Although we develop \TriCG and \TriMR for SQD systems, both can be generalized to handle any linear system of the form
  \begin{equation}
    \label{eq:spd_snd_sqd}
    \begin{bmatrix}
      \pm M   &     A \\
          A^T & \pm N
    \end{bmatrix}
    \begin{bmatrix}
    x \\
    y
    \end{bmatrix}
    =
    \begin{bmatrix}
    b \\
    c
    \end{bmatrix},
  \end{equation}
  where \(M\) and \(N\) are symmetric positive definite.
  Symmetric definite systems can always be written in the form~\eqref{eq:spd_snd_sqd}.
  For instance, one could solve any definite system by first partitioning any symmetric permutation of it as~\eqref{eq:spd_snd_sqd} and applying \TriCG or \TriMR to the resulting \(A\).
  There are multiple ways to perform such partitioning, and it is not clear whether such strategy might lead to improved solution processes for SPD systems.
  With minor modifications, \TriMR also supports the saddle-point system~\eqref{eq:saddle_point}.
  In this case, $N$ can be replaced by any SPD matrix to define an elliptic norm in \Cref{alg:ssy}, which extends the possibilities for preconditioning the linear system. For instance, \TriMR could use the preconditioner $\blkdiag(M, S)$ \citep{murphy-golub-wathen-2000} for Stokes problems discretized by LBB-stable finite element pairs where $S$ is the Schur complement $A\T M^{-1} A$ or an approximation to it.
  Our implementations of \TriCG and \TriMR take into account all these extensions, are applicable in any floating-point system supported by Julia, and run on GPUs.

  \TriCG and \TriMR perform substantially better than \SYMMLQ and \MINRES in our experiments and often terminate in about half as many iterations on a residual-based stopping condition.
  % This behavior is predictable because block Krylov spaces enlarge the search space.
  Nevertheless more extensive testing is required to properly assess their performance in practice.
  Although theoretically equivalent, \TriCG and \TriMR appear to preserve orthogonality of the Krylov basis better than their counterparts \BLCG and \BLMINRES.
  Additional numerical illustrations are available in \cref{sec:supplements}.

  Based upon \Cref{alg:ssy}, it is possible to develop a third method in the spirit of \SYMMLQ that we could name \TriLQ.
  The \TriLQ subproblem selects \(z_k\) in~\eqref{eq:xy} as the solution of the minimum-norm subproblem
  \[
    \minimize{z_k \in \R^{2k}} \ \|z_k\|
    \quad
    \st \ S_{k-1, k} z_k = \beta_1 e_1 + \gamma_1 e_2,
  \]
  where \(S_{k-1,k}\) is the leading \((2k-2)\)\(\times\)\((2k)\) submatrix of \(S_{k+1,k}\) in~\eqref{eq:sqd-block-lanczos}.
  The subproblem can be solved via the LQ factorization of \(S_{k-1,k}\).
  Much of \TriLQ would be similar to block-\SYMMLQ: iterates are updated along orthogonal directions, the \(H\)-norm of the iterates increases monotonically, and the \(H\)-norm of the error decreases monotonically.
  At each iteration, \TriLQ allows the user to transfer to the \TriCG iterate.
  Because the \TriCG iterate always exists for~\eqref{eq:sqd}, \TriLQ might not have have any advantage in practice, other than completing the family of numerical methods based on \Cref{alg:ssy}.
  However, the \TriLQ iterate remains well defined for the saddle-point system~\eqref{eq:saddle_point}, whereas \TriCG may break down in that case.

  The strong connection between \TriCG and block-\CGM with blocks of size \(2\) suggests that \TriCG might also be useful to approximate eigenvalues.
  We leave the investigation of such extensions to future work.

  \subsection*{Acknowledgements}

  We sincerely thank Michael A. Saunders and two anonymous referees for numerous suggestions that improved the content and presentation of the present research.

  \small
  \bibliographystyle{abbrvnat}
  \bibliography{abbrv,tricg}
  \normalsize

  \clearpage
  \appendix

  \section{Additional numerical results}%
  \label{sec:supplements}

\begin{figure}[ht]
    \centering
    \includetikzgraphics[width=0.49\textwidth]{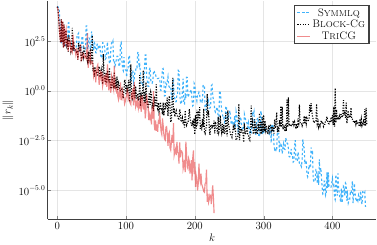}
    \hfill
    \includetikzgraphics[width=0.49\textwidth]{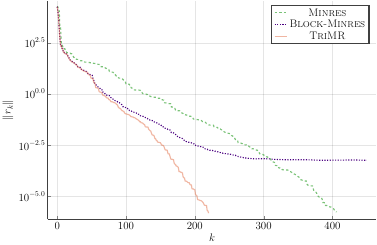}
    \label{fig:lp_d6cube}
    \caption{Residual history of \SYMMLQ, \BLCG, \TriCG, \MINRES, \BLMINRES and \TriMR on lp\_d6cube, where $A$ is $415 \times 6184$ with $37704$ nonzeros.}
  \end{figure}

\begin{figure}[ht]
    \centering
    \includetikzgraphics[width=0.49\textwidth]{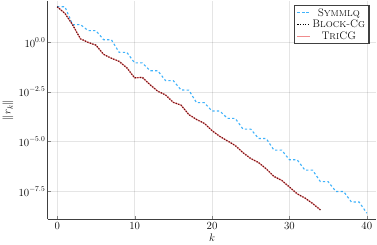}
    \hfill
    \includetikzgraphics[width=0.49\textwidth]{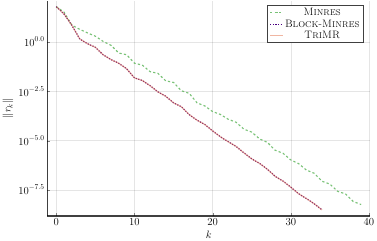}
    \label{fig:well1033}
    \caption{Residual history of \SYMMLQ, \BLCG, \TriCG, \MINRES, \BLMINRES and \TriMR on well1033, where $A$ is $1033 \times 320$ with $4732$ nonzeros.}
  \end{figure}

\begin{figure}[ht]
    \centering
    \includetikzgraphics[width=0.49\textwidth]{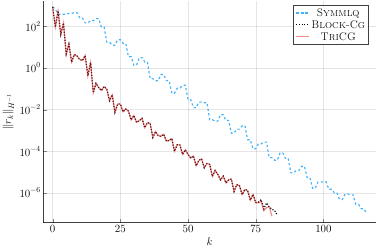}
    \hfill
    \includetikzgraphics[width=0.49\textwidth]{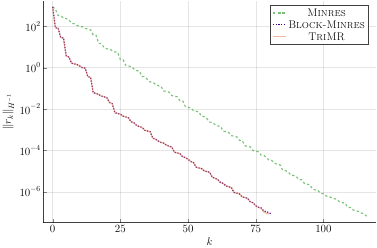}
    \label{fig:colliding_flow1}
    \caption{Residual history of \SYMMLQ, \BLCG, \TriCG, \MINRES, \BLMINRES and \TriMR iterates on colliding\_flow problem with \(Q_1-P_0\) discretization.
    The discretized linear system has size $12546 \times 12546$ with $146241$ nonzeros.}
  \end{figure}

\begin{figure}[ht]
    \centering
    \includetikzgraphics[width=0.49\textwidth]{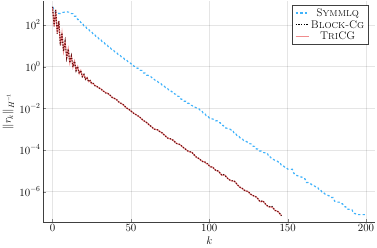}
    \hfill
    \includetikzgraphics[width=0.49\textwidth]{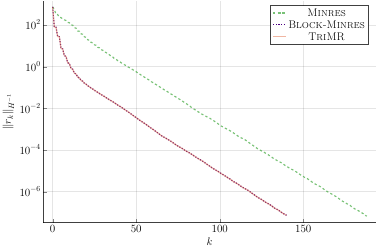}
    \label{fig:colliding_flow2}
    \caption{Residual history of \SYMMLQ, \BLCG, \TriCG, \MINRES, \BLMINRES and \TriMR iterates on colliding\_flow problem with \(Q_1-Q_1\) discretization.
    The discretized linear system has size $12675 \times 12675$ with $239873$ nonzeros.}
  \end{figure}

\begin{figure}[ht]
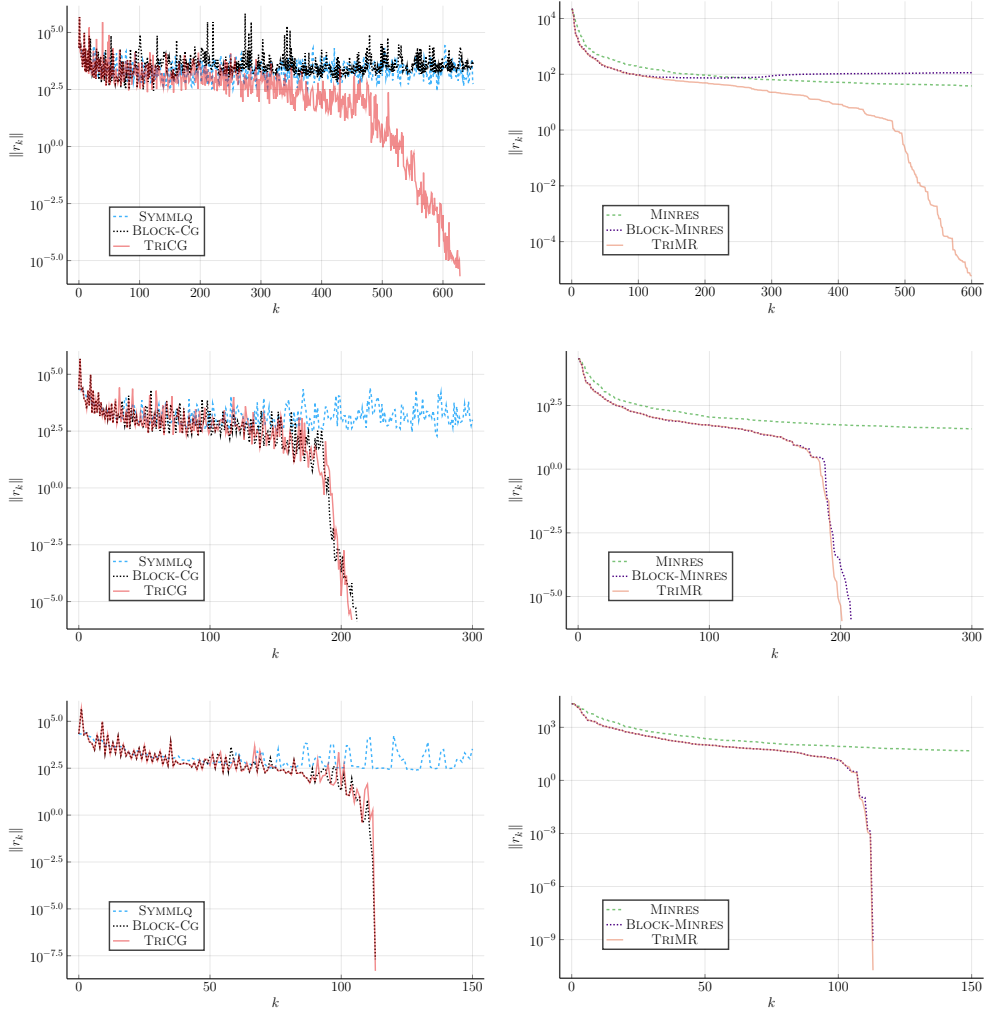

    \centering
    \includetikzgraphics[width=0.49\textwidth]{mp/lpi_klein3_tricg_blcg_symmlq_Float64}
    \hfill
    \includetikzgraphics[width=0.49\textwidth]{mp/lpi_klein3_trimr_blminres_minres_Float64}
    \vfill
    \includetikzgraphics[width=0.49\textwidth]{mp/lpi_klein3_tricg_blcg_symmlq_Float128}
    \hfill
    \includetikzgraphics[width=0.49\textwidth]{mp/lpi_klein3_trimr_blminres_minres_Float128}
    \vfill
    \includetikzgraphics[width=0.49\textwidth]{mp/lpi_klein3_tricg_blcg_symmlq_Float256}
    \hfill
    \includetikzgraphics[width=0.49\textwidth]{mp/lpi_klein3_trimr_blminres_minres_Float256}
    \label{fig:lpi_klein3_mp}
    \caption{Residual history of \SYMMLQ, \BLCG, \TriCG, \MINRES, \BLMINRES and \TriMR iterates on lpi\_klein3 in double (top), quadruple (center) and octuple precision (bottom) where $A$ is $994 \times 1082$ with $13101$ nonzeros.}
  \end{figure}

  This appendix contains further numerical comparisons between \TriCG, \TriMR, \BLCG, \BLMINRES, \SYMMLQ and \MINRES.
  \Cref{fig:lp_d6cube} reports residual histories in double precision on another underdetermined system from optimization, where \BLCG and \BLMINRES do not converge, presumably due to excessive loss of orthogonality.
  \Cref{fig:well1033} corresponds to a well-conditioned overdetermined system from a least-squares application, where the residuals of the block methods nearly coincide.
  \Cref{fig:colliding_flow1,fig:colliding_flow2} are Stokes systems.
  Finally, \Cref{fig:lpi_klein3_mp} is a rather dramatic example of an underdetermined system from optimization where only \TriCG and \TriMR converge in double precision.
  As the accuracy increases, \BLCG and \BLMINRES converge and nearly coincide with \TriCG and \TriMR.
  Moreover, as the number of digits doubles, the number of iterations to converge is roughly halved.

  \begin{landscape}
  \section{\TriMR details}%
  \label{sec:appendix}

  \begin{equation}
  \label{eq:Q}
  Q_{2k-1,2k+2}
  :=
  \bmat{
  1    &         &                     &      \\
       & c_{4,k} & \phantom{-} s_{4,k} &      \\
       & s_{4,k} &          -  c_{4,k} &      \\
       &         &                     & 1
  }
  \bmat{
  1    &         &      &      \\
       & c_{3,k} &      & \phantom{-} s_{3,k} \\
       &         & 1    &                     \\
       & s_{3,k} &      &          -  c_{3,k}
  }
  \bmat{
  c_{2,k} & \phantom{-} s_{2,k} &      &      \\
  s_{2,k} &          -  c_{2,k} &      &      \\
          &                     & 1    &      \\
          &                     &      & 1
  }
  \bmat{
  1    &         &      &      \\
       & c_{1,k} &      & \phantom{-} s_{1,k} \\
       &         & 1    &                     \\
       & s_{1,k} &      &          -  c_{1,k}
  }
  \end{equation}

  {\small
  \[
    \kbordermatrix{
         & 2k-1 & 2k      & 2k+1 & 2k+2 \\
    2k-1 & 1    &         &      &      \\
    2k   &      & c_{1,k} &      & \phantom{-} s_{1,k} \\
    2k+1 &      &         & 1    &                     \\
    2k+2 &      & s_{1,k} &      &          -  c_{1,k}
    }
    \!\!\!\!
    \kbordermatrix{
    & 2k-1                & 2k                  &              & 2k+1               & 2k+2                 & 2k+3        & 2k+4         \\
    & \bar{\delta}_{2k-1} & \bar{\sigma}_{2k-1} & \! \vrule \! & \bar{\eta}_{2k-1}  & \bar{\lambda}_{2k-1} & 0           & 0            \\
    & \bar{\theta}_{k}    & \bar{\delta}_{2k}   & \! \vrule \! &  \bar{\sigma}_{2k} & 0                    & 0           & 0            \\
    & 0                   & \gamma_{k+1}        & \! \vrule \! &  1                 & \alpha_{k+1}         & 0           & \gamma_{k+2} \\
    & \beta_{k+1}         & 0                   & \! \vrule \! &  \alpha_{k+1}      & -1                   & \beta_{k+2} & 0
    }
    =
    \!\!\!\!
    \kbordermatrix{
    & 2k-1                & 2k                        &              & 2k+1                     & 2k+2                      & 2k+3                      & 2k+4         \\
    & \bar{\delta}_{2k-1} & \bar{\sigma}_{2k-1}       & \! \vrule \! & \bar{\eta}_{2k-1}        & \bar{\lambda}_{2k-1}      & 0                         & 0            \\
    & \theta_{k}          & \widetilde{\delta}_{2k}   & \! \vrule \! & \widetilde{\sigma}_{2k}  & \widetilde{\eta}_{2k}     & \widetilde{\lambda}_{2k}  & 0            \\
    & 0                   & \gamma_{k+1}              & \! \vrule \! &  1                       & \alpha_{k+1}              & 0                         & \gamma_{k+2} \\
    & 0                   & g_k                       & \! \vrule \! & \widetilde{\theta}_{k+1} & \widetilde{\delta}_{2k+2} & \widetilde{\sigma}_{2k+2} & 0
    }
  \]
  }

  \begin{subequations}
    \label{eq:ref1}
    \begin{equation}
      \theta_{k} = \sqrt{\bar{\theta}_{k}^2 + \beta_{k+1}^2},\quad
      c_{1,k} = \theta_{k} / \bar{\theta}_{k},\quad
      s_{1,k} = \beta_{k+1} / \bar{\theta}_k
    \end{equation}
    \begin{align}
        \widetilde{\delta}_{2k}   &= c_{1,k} \bar{\delta}_{2k}                       ,
      & \widetilde{\sigma}_{2k}   &= c_{1,k} \bar{\sigma}_{2k} + s_{1,k} \alpha_{k+1},
      & \widetilde{\eta}_{2k}     &=                           - s_{1,k}             ,
      & \widetilde{\lambda}_{2k}  &=                 \phantom{-} s_{1,k} \beta_{k+2} , \\
        g_k                       &= s_{1,k} \bar{\delta}_{2k}                       ,
      & \widetilde{\theta}_{k+1}  &= s_{1,k} \bar{\sigma}_{2k} - c_{1,k} \alpha_{k+1},
      & \widetilde{\delta}_{2k+2} &=                 \phantom{-} c_{1,k}             ,
      & \widetilde{\sigma}_{2k+2} &=                           - c_{1,k} \beta_{k+2} .
    \end{align}
  \end{subequations}

  {\small
  \[
    \kbordermatrix{
         & 2k-1    & 2k                  & 2k+1 & 2k+2 \\
    2k-1 & c_{2,k} & \phantom{-} s_{2,k} &      &      \\
    2k   & s_{2,k} &          -  c_{2,k} &      &      \\
    2k+1 &         &                     & 1    &      \\
    2k+2 &         &                     &      & 1
    }
    \!\!\!\!
    \kbordermatrix{
    & 2k-1                & 2k                        &              & 2k+1                     & 2k+2                      & 2k+3                      & 2k+4         \\
    & \bar{\delta}_{2k-1} & \bar{\sigma}_{2k-1}       & \! \vrule \! & \bar{\eta}_{2k-1}        & \bar{\lambda}_{2k-1}      & 0                         & 0            \\
    & \theta_{k}          & \widetilde{\delta}_{2k}   & \! \vrule \! & \widetilde{\sigma}_{2k}  & \widetilde{\eta}_{2k}     & \widetilde{\lambda}_{2k}  & 0            \\
    & 0                   & \gamma_{k+1}              & \! \vrule \! &  1                       & \alpha_{k+1}              & 0                         & \gamma_{k+2} \\
    & 0                   & g_k                       & \! \vrule \! & \widetilde{\theta}_{k+1} & \widetilde{\delta}_{2k+2} & \widetilde{\sigma}_{2k+2} & 0
    }
    =
    \!\!\!\!
    \kbordermatrix{
    & 2k-1          & 2k                    &              & 2k+1                     & 2k+2                      & 2k+3                      & 2k+4         \\
    & \delta_{2k-1} & \sigma_{2k-1}         & \! \vrule \! & \eta_{2k-1}              & \lambda_{2k-1}            & \mu_{2k-1}                & 0            \\
    & 0             & \widehat{\delta}_{2k} & \! \vrule \! & \widehat{\sigma}_{2k}    & \widehat{\eta}_{2k}       & \widehat{\lambda}_{2k}    & 0            \\
    & 0             & \gamma_{k+1}          & \! \vrule \! &  1                       & \alpha_{k+1}              & 0                         & \gamma_{k+2} \\
    & 0             & g_k                   & \! \vrule \! & \widetilde{\theta}_{k+1} & \widetilde{\delta}_{2k+2} & \widetilde{\sigma}_{2k+2} & 0
    }
  \]
  }

  \begin{subequations}
    \label{eq:ref2}
    \begin{equation}
      \delta_{2k-1} = \sqrt{\bar{\delta}_{2k-1}^2 + \theta_{k}^2},\quad c_{2,k} = \bar{\delta}_{2k-1} / \delta_{2k-1},\quad s_{2,k} = \theta_{k} / \delta_{2k-1}
    \end{equation}
    \begin{align}
        \sigma_{2k-1}            &= c_{2,k} \bar{\sigma}_{2k-1}  + s_{2,k} \widetilde{\delta}_{2k} ,
      & \eta_{2k-1}              &= c_{2,k} \bar{\eta}_{2k-1}    + s_{2,k} \widetilde{\sigma}_{2k} ,
      & \lambda_{2k-1}           &= c_{2,k} \bar{\lambda}_{2k-1} + s_{2,k} \widetilde{\eta}_{2k}   ,
      & \mu_{2k-1}               &=                    \phantom{-} s_{2,k} \widetilde{\lambda}_{2k}, \\
        \widehat{\delta}_{2k}    &= s_{2,k} \bar{\sigma}_{2k-1}  - c_{2,k} \widetilde{\delta}_{2k} ,
      & \widehat{\sigma}_{2k}    &= s_{2,k} \bar{\eta}_{2k-1}    - c_{2,k} \widetilde{\sigma}_{2k} ,
      & \widetilde{\eta}_{2k}    &= s_{2,k} \bar{\lambda}_{2k-1} - c_{2,k} \widetilde{\eta}_{2k}   ,
      & \widetilde{\lambda}_{2k} &=                              - c_{2,k} \widetilde{\lambda}_{2k}.
    \end{align}
  \end{subequations}

  {\small
  \[
    \kbordermatrix{
         & 2k-1 & 2k      & 2k+1 & 2k+2 \\
    2k-1 & 1    &         &      &      \\
    2k   &      & c_{3,k} &      & \phantom{-} s_{3,k} \\
    2k+1 &      &         & 1    &                     \\
    2k+2 &      & s_{3,k} &      &          -  c_{3,k}
    }
    \!\!\!\!
    \kbordermatrix{
    & 2k-1          & 2k                    &              & 2k+1                     & 2k+2                      & 2k+3                      & 2k+4         \\
    & \delta_{2k-1} & \sigma_{2k-1}         & \! \vrule \! & \eta_{2k-1}              & \lambda_{2k-1}            & \mu_{2k-1}                & 0            \\
    & 0             & \widehat{\delta}_{2k} & \! \vrule \! & \widehat{\sigma}_{2k}    & \widehat{\eta}_{2k}       & \widehat{\lambda}_{2k}    & 0            \\
    & 0             & \gamma_{k+1}          & \! \vrule \! &  1                       & \alpha_{k+1}              & 0                         & \gamma_{k+2} \\
    & 0             & g_k                   & \! \vrule \! & \widetilde{\theta}_{k+1} & \widetilde{\delta}_{2k+2} & \widetilde{\sigma}_{2k+2} & 0
    }
    =
    \!\!\!\!
    \kbordermatrix{
    & 2k-1          & 2k                     &              & 2k+1                   & 2k+2                 & 2k+3                    & 2k+4         \\
    & \delta_{2k-1} & \sigma_{2k-1}          & \! \vrule \! & \eta_{2k-1}            & \lambda_{2k-1}       & \mu_{2k-1}              & 0            \\
    & 0             & \mathring{\delta}_{2k} & \! \vrule \! & \mathring{\sigma}_{2k} & \mathring{\eta}_{2k} & \mathring{\lambda}_{2k} & 0            \\
    & 0             & \gamma_{k+1}           & \! \vrule \! &  1                     & \alpha_{k+1}         & 0                       & \gamma_{k+2} \\
    & 0             & 0                      & \! \vrule \! & \bar{\theta}_{k+1}     & \bar{\delta}_{2k+2}  & \bar{\sigma}_{2k+2}     & 0
    }
  \]
  }

  \begin{subequations}
    \label{eq:ref3}
    \begin{equation}
      \mathring{\delta}_{2k} = \sqrt{\widehat{\delta}_{2k}^2 + g_k^2},\quad c_{3,k} = \widehat{\delta}_{2k} / \mathring{\delta}_{2k},\quad s_{3,k} = g_k / \mathring{\delta}_{2k}
    \end{equation}
    \begin{align}
        \mathring{\sigma}_{2k}  &= c_{3,k} \widehat{\sigma}_{2k}  + s_{3,k} \widetilde{\theta}_{k+1} ,
      & \mathring{\eta}_{2k}    &= c_{3,k} \widehat{\eta}_{2k}    + s_{3,k} \widetilde{\delta}_{2k+2},
      & \mathring{\lambda}_{2k} &= c_{3,k} \widehat{\lambda}_{2k} + s_{3,k} \widetilde{\sigma}_{2k+2}, \\
        \bar{\theta}_{k+1}      &= s_{3,k} \widehat{\sigma}_{2k}  - c_{3,k} \widetilde{\theta}_{k+1} ,
      & \bar{\delta}_{2k+2}     &= s_{3,k} \widehat{\eta}_{2k}    - c_{3,k} \widetilde{\delta}_{2k+2},
      & \bar{\sigma}_{2k+2}     &= s_{3,k} \widehat{\lambda}_{2k} - c_{3,k} \widetilde{\sigma}_{2k+2}.
    \end{align}
  \end{subequations}

  {\small
  \[
    \kbordermatrix{
         & 2k-1 & 2k      & 2k+1                & 2k+2 \\
    2k-1 & 1    &         &                     &      \\
    2k   &      & c_{4,k} & \phantom{-} s_{4,k} &      \\
    2k+1 &      & s_{4,k} &          -  c_{4,k} &      \\
    2k+2 &      &         &                     & 1
    }
    \!\!\!\!
    \kbordermatrix{
    & 2k-1          & 2k                     &              & 2k+1                     & 2k+2                      & 2k+3                      & 2k+4         \\
    & \delta_{2k-1} & \sigma_{2k-1}          & \! \vrule \! & \eta_{2k-1}              & \lambda_{2k-1}            & \mu_{2k-1}                & 0            \\
    & 0             & \mathring{\delta}_{2k} & \! \vrule \! & \mathring{\sigma}_{2k}   & \mathring{\eta}_{2k}      & \mathring{\lambda}_{2k}   & 0            \\
    & 0             & \gamma_{k+1}           & \! \vrule \! &  1                       & \alpha_{k+1}              & 0                         & \gamma_{k+2} \\
    & 0             & 0                      & \! \vrule \! & \widetilde{\theta}_{k+1} & \widetilde{\delta}_{2k+2} & \widetilde{\sigma}_{2k+2} & 0
    }
    =
    \!\!\!\!
    \kbordermatrix{
    & 2k-1          & 2k            &              & 2k+1                & 2k+2                & 2k+3                & 2k+4                 \\
    & \delta_{2k-1} & \sigma_{2k-1} & \! \vrule \! & \eta_{2k-1}         & \lambda_{2k-1}      & \mu_{2k-1}          & 0                    \\
    & 0             & \delta_{2k}   & \! \vrule \! & \sigma_{2k}         & \eta_{2k}           & \lambda_{2k}        & \mu_{2k}             \\
    & 0             & 0             & \! \vrule \! & \bar{\delta}_{2k+1} & \bar{\sigma}_{2k+1} & \bar{\eta}_{2k+1}   & \bar{\lambda}_{2k+1} \\
    & 0             & 0             & \! \vrule \! & \bar{\theta}_{k+1}  & \bar{\delta}_{2k+2} & \bar{\sigma}_{2k+2} & 0
    }
  \]
  }

  \begin{subequations}
    \label{eq:ref4}
    \begin{equation}
      \delta_{2k} = \sqrt{\mathring{\delta}_{2k}^2 + \gamma_{k+1}^2},\quad c_{4,k} = \mathring{\delta}_{2k}^2 / \delta_{2k},\quad s_{4,k} = \gamma_{k+1} / \delta_{2k}
    \end{equation}
    \begin{align}
        \sigma_{2k}          &= c_{4,k} \mathring{\sigma}_{2k}  + s_{4,k}             ,
      & \eta_{2k}            &= c_{4,k} \mathring{\eta}_{2k}    + s_{4,k} \alpha_{k+1},
      & \lambda_{2k}         &= c_{4,k} \mathring{\lambda}_{2k}                       ,
      & \mu_{2k}             &=                       \phantom{-} s_{4,k} \gamma_{k+2}, \\
        \bar{\delta}_{2k+1}  &= s_{4,k} \mathring{\sigma}_{2k}  - c_{4,k}             ,
      & \bar{\sigma}_{2k+1}  &= s_{4,k} \mathring{\eta}_{2k}    - c_{4,k} \alpha_{k+1},
      & \bar{\eta}_{2k+1}    &= s_{4,k} \mathring{\lambda}_{2k}                       ,
      & \bar{\lambda}_{2k+1} &=                                 - c_{4,k} \gamma_{k+2}.
    \end{align}
  \end{subequations}

  \[
    \bmat{
    1    &         &      &      \\
         & c_{1,k} &      & \phantom{-} s_{1,k} \\
         &         & 1    &                     \\
         & s_{1,k} &      &          -  c_{1,k}
    }
    \bmat{\bar{\pi}_{2k-1} \\ \bar{\pi}_{2k} \\ 0 \\ 0}
    =
    \bmat{\bar{\pi}_{2k-1} \\ \widetilde{\pi}_{2k} \\ 0 \\ \widetilde{\pi}_{2k+2}}
    %\bmat{\bar{\pi}_{2k-1} \\ c_{1,k} \bar{\pi}_{2k} \\ 0 \\ s_{1,k} \bar{\pi}_{2k}}
    \qquad \qquad \qquad
    \bmat{
    c_{2,k} & \phantom{-} s_{2,k} &      &      \\
    s_{2,k} &          -  c_{2,k} &      &      \\
            &                     & 1    &      \\
            &                     &      & 1
    }
    \bmat{\bar{\pi}_{2k-1} \\ \widetilde{\pi}_{2k} \\ 0 \\ \widetilde{\pi}_{2k+2}}
    %\bmat{\bar{\pi}_{2k-1} \\ c_{1,k} \bar{\pi}_{2k} \\ 0 \\ s_{1,k} \bar{\pi}_{2k}}
    =
    \bmat{\pi_{2k-1} \\ \widehat{\pi}_{2k} \\ 0 \\ \widetilde{\pi}_{2k+2}}
    %\bmat{c_{2,k} \bar{\pi}_{2k-1} + s_{2,k} c_{1,k} \bar{\pi}_{2k} \\ s_{2,k} \bar{\pi}_{2k-1} - c_{2,k} c_{1,k} \bar{\pi}_{2k} \\ 0 \\ s_{1,k} \bar{\pi}_{2k}}
  \]

  \begin{equation}
    \label{eq:ref12}
    \widetilde{\pi}_{2k} = c_{1,k} \bar{\pi}_{2k},\quad \widetilde{\pi}_{2k+2} = s_{1,k} \bar{\pi}_{2k}, \quad \pi_{2k-1} = c_{2,k} \bar{\pi}_{2k-1} + s_{2,k} \widetilde{\pi}_{2k},\quad \widehat{\pi}_{2k} = s_{2,k} \bar{\pi}_{2k-1} - c_{2,k} \widetilde{\pi}_{2k}
  \end{equation}

  \[
    \bmat{
    1    &         &      &      \\
         & c_{3,k} &      & \phantom{-} s_{3,k} \\
         &         & 1    &                     \\
         & s_{3,k} &      &          -  c_{3,k}
    }
    \bmat{\pi_{2k-1} \\ \widehat{\pi}_{2k} \\ 0 \\ \widetilde{\pi}_{2k+2}}
    %\bmat{c_{2,k} \bar{\pi}_{2k-1} + s_{2,k} c_{1,k} \bar{\pi}_{2k} \\ s_{2,k} \bar{\pi}_{2k-1} - c_{2,k} c_{1,k} \bar{\pi}_{2k} \\ 0 \\ s_{1,k} \bar{\pi}_{2k}}
    =
    \bmat{\pi_{2k-1} \\ \mathring{\pi}_{2k} \\ 0 \\ \bar{\pi}_{2k+2}}
    %\bmat{c_{2,k} \bar{\pi}_{2k-1} + s_{2,k} c_{1,k} \bar{\pi}_{2k} \\ c_{3,k} s_{2,k} \bar{\pi}_{2k-1} + ( s_{3,k} s_{1,k} -  c_{3,k} c_{2,k} c_{1,k}) \bar{\pi}_{2k} \\ 0 \\ s_{3,k} s_{2,k} \bar{\pi}_{2k-1} - (c_{3,k} s_{1,k} + s_{3,k} c_{2,k} c_{1,k}) \bar{\pi}_{2k}}
    \qquad \qquad \qquad
    \bmat{
    1    &         &                     &      \\
         & c_{4,k} & \phantom{-} s_{4,k} &      \\
         & s_{4,k} &          -  c_{4,k} &      \\
         &         &                     & 1
    }
    \bmat{\pi_{2k-1} \\ \mathring{\pi}_{2k} \\ 0 \\ \bar{\pi}_{2k+2}}
    %\bmat{c_{2,k} \bar{\pi}_{2k-1} + s_{2,k} c_{1,k} \bar{\pi}_{2k} \\ c_{3,k} s_{2,k} \bar{\pi}_{2k-1} + ( s_{3,k} s_{1,k} -  c_{3,k} c_{2,k} c_{1,k}) \bar{\pi}_{2k} \\ 0 \\ s_{3,k} s_{2,k} \bar{\pi}_{2k-1} - (c_{3,k} s_{1,k} + s_{3,k} c_{2,k} c_{1,k}) \bar{\pi}_{2k}}
    =
    \bmat{\pi_{2k-1} \\ \pi_{2k} \\ \bar{\pi}_{2k+1} \\ \bar{\pi}_{2k+2}}
    %\bmat{c_{2,k} \bar{\pi}_{2k-1} + s_{2,k} c_{1,k} \bar{\pi}_{2k} \\ c_{4,k} c_{3,k} s_{2,k} \bar{\pi}_{2k-1} + c_{4,k} ( s_{3,k} s_{1,k} -  c_{3,k} c_{2,k} c_{1,k}) \bar{\pi}_{2k} \\ s_{4,k} c_{3,k} s_{2,k} \bar{\pi}_{2k-1} + s_{4,k} ( s_{3,k} s_{1,k} -  c_{3,k} c_{2,k} c_{1,k}) \bar{\pi}_{2k} \\ s_{3,k} s_{2,k} \bar{\pi}_{2k-1} - (c_{3,k} s_{1,k} + s_{3,k} c_{2,k} c_{1,k}) \bar{\pi}_{2k}}
  \]

  \begin{equation}
    \label{eq:ref34}
    \mathring{\pi}_{2k} = c_{3,k} \widehat{\pi}_{2k} + s_{3,k} \widetilde{\pi}_{2k+2},\quad \bar{\pi}_{2k+2} = s_{3,k} \widehat{\pi}_{2k} - c_{3,k} \widetilde{\pi}_{2k+2}, \quad
  \pi_{2k} = c_{4,k} \mathring{\pi}_{2k},\quad \bar{\pi}_{2k+1} = s_{4,k} \mathring{\pi}_{2k}
  \end{equation}

  % \[\pi_{2k-1} = c_{2,k} \bar{\pi}_{2k-1} + s_{2,k} c_{1,k} \bar{\pi}_{2k}\]
  % \[\pi_{2k} = c_{4,k} c_{3,k} s_{2,k} \bar{\pi}_{2k-1} + c_{4,k} ( s_{3,k} s_{1,k} -  c_{3,k} c_{2,k} c_{1,k}) \bar{\pi}_{2k}\]
  % \[\bar{\pi}_{2k+1} = s_{4,k} c_{3,k} s_{2,k} \bar{\pi}_{2k-1} + s_{4,k} ( s_{3,k} s_{1,k} -  c_{3,k} c_{2,k} c_{1,k}) \bar{\pi}_{2k}\]
  % \[\bar{\pi}_{2k+2} = s_{3,k} s_{2,k} \bar{\pi}_{2k-1} - (c_{3,k} s_{1,k} + s_{3,k} c_{2,k} c_{1,k}) \bar{\pi}_{2k}\]

  Scalars decorated by a hat, a tilde or a ring are updated at the current iteration.
  Scalars decorated by a bar will be updated at the next iteration.
  Scalars without any decoration have been updated to their final value.
  \end{landscape}

\end{document}